\documentclass[12pt]{article}

\usepackage{amssymb,amsthm,amsmath,amssymb}
\usepackage{graphicx}
\usepackage[colorlinks=true,citecolor=black,linkcolor=black,urlcolor=blue]{hyperref}

\usepackage{xcolor}

\usepackage{pstricks,pst-node,enumitem}

\newtheorem{thm}{Theorem}
\newtheorem{proposition}[thm]{Proposition}

\newtheorem{corollary}[thm]{Corollary}

\newtheorem{remark}[thm]{Remark}

\usepackage{fullpage}

\DeclareMathOperator{\hv}{\mathrm{hv}}
\DeclareMathOperator{\lv}{\mathrm{lv}}
\DeclareMathOperator{\V}{\mathrm{v}}
\newcommand{\1}{\mathbf{1}}

\title{Chains of binary paths and shifted tableaux}

\author{K. Manes, I. Tasoulas, A. Sapounakis, P. Tsikouras\\
\small Department of Informatics, University of Piraeus\\ 
\small Piraeus, GREECE \\
\small\tt \{kmanes, jtas, arissap, pgtsik\}@unipi.gr}

\begin{document}

\maketitle

\begin{abstract}
In this paper, a natural bijection between multichains of binary paths and shifted tableaux is presented, and it is used for the enumeration of the chains with maximum length from a given path $P$ to the maximum path $\1_{|P|}$.
By mapping chains to shifted tableaux, the main formulas given in a recent paper by the authors for the enumeration of the $P - \1_{|P|}$ chains having only small intervals and minimum length are proved, using some new bijections on shifted tableaux.

\end{abstract}

\section{Introduction}\label{section:intro}

Let $\mathcal{P}_n$ be the set of all (binary) paths $P$ of length $|P| = n$, i.e., lattice paths $P = p_1 p_2 \cdots p_n$ where each {\em step} $p_i$, $i \in [n]$, is either an {\em upstep} $u = (1,1)$ or a {\em downstep} $d = (1,-1)$ and connects two consecutive points of the path. We denote by $|P|_u$ (resp. $|P|_d$) the number of upsteps (resp. downsteps) of $P$. An {\em ascent} (resp. {\em descent}) of $P$ is a maximal sequence of $u$'s (resp. $d$'s) in $P$. The last point of an ascent (resp. descent) is called {\em peak} (resp. {\em valley}) of the path. Clearly, every peak (resp. valley) corresponds to either an occurrence of $ud$ (resp. $du$), or an occurrence of $u$ (resp. $d$) at the end of the path. It is convenient to consider that the starting point of a path is the origin of a pair of axes. The {\em height} of a point of $P$ is its $y$-coordinate. We denote by $\lv(P)$ (resp. $\hv(P)$) the height of the lowest (resp. highest) valley of $P$. A valley of $P$ with height $\lv(P)$ is called {\em low} valley of $P$. We set $\mathcal{P} = \bigcup_{n \ge 0} \mathcal{P}_n$, where $\mathcal{P}_0$ consists of only the empty path $\varepsilon$ (the path which has no steps). 

A {\em Dyck path} is a path that starts and ends at the same height and lies weakly above this height. In this paper, all Dyck paths are denoted by lowercase letters. The set of Dyck paths of length $2n$ is denoted by $\mathcal{D}_n$, and we set $\mathcal{D} = \bigcup_{n \ge 0} \mathcal{D}_n$, where $\mathcal{D}_0 = \{\varepsilon\}$. 

A path which is a prefix (resp. suffix) of a Dyck path, is called {\em Dyck prefix} (resp. {\em Dyck suffix}). Every non-initial point of a Dyck prefix having height zero is called {\em return} point. A Dyck path with only one return point is called {\em prime}. Every non-empty Dyck path $a$ can be uniquely decomposed as a product of prime Dyck paths, i.e., $a = u a_1 d u a_2 d \cdots u a_k d$, where $a_i \in \mathcal{D}$, $i \in [k]$. Furthermore, every Dyck prefix (resp. Dyck suffix) $P$ can be uniquely decomposed in the form $P = a_0 u a_1 \cdots u a_k$ (resp. $P = a_0 d a_1 \cdots d a_k$), where $a_i \in \mathcal{D}$, $i \in [0,k]$, $k \ge 0$. 

A natural (partial) ordering on $\mathcal{P}_n$ is defined via the geometric representation of the paths $P, Q \in \mathcal{P}_n$, where $P \le Q$ whenever $P$ lies (weakly) below $Q$. Obviously, a path $Q$ covers a path $P$ whenever $Q$ is obtained from $P$ by turning exactly one valley of $P$ into a peak. 
This ordering is better understood by considering the following alternative encoding of binary paths: Every $P \in \mathcal{P}_n$ can be described uniquely by the sequence $(h_i(P))_{i \in [n]}$ of the heights of its points, so that $P \le Q$ iff $h_i(P) \le h_i(Q)$, $i \in [n]$. Then, the join and meet of $P$, $Q$ are given by 
\[ h_i(P \vee Q) = \max\{ h_i(P), h_i(Q) \} \textrm{ and } h_i(P \wedge Q) = \min\{ h_i(P), h_i(Q) \}. \]
It is well-known that the poset $(\mathcal{P}_n, \le)$ is a finite, self-dual, distributive lattice with minimum and maximum elements the paths $\mathbf{0}_n = d^n = \underbrace{d d \cdots d}_{n \ \rm{times}}$ and $\mathbf{1}_n = u^n = \underbrace{u u \cdots u}_{n \ \rm{times}}$ respectively.
We note that the length of every saturated $P - Q$ chain, where $P = p_1 p_2 \cdots p_n$, $Q = q_1 q_2 \cdots q_n$, is equal to 
\[ l(P,Q) = \frac{1}{2} \sum\limits_{i=1}^n (h_i(Q) - h_i(P)) = \sum\limits_{i=1}^n (n-i+1) \cdot ([q_i = u] - [p_i = u]), \]
where $[S]$ is the Iverson binary notation, i.e., for every proposition $S$, $[S] = 1$ if $S$ is true, and $0$ if $S$ is false. 

The above lattice appears in the literature in various equivalent forms (e.g., binary words \cite[p. 92]{KnuthV4A}, subsets of $[n]$ \cite{Lindstrom1970}, permutations of $[n]$ \cite[p. 402]{Stanley2011},  partitions of $n$ into distinct parts \cite{Stanley1991}, threshold graphs \cite{MerrisRoby2005}). Sapounakis et al \cite{STT2006} studied its sublattice of Dyck paths. Recently, a bijection between comparable pairs of paths of this lattice and Dyck prefixes of odd length has been presented \cite{MSTT2018}.

An interval $[x,y]$ in a poset is called {\em small} if $y$ is the
 join of some elements covering $x$. It is easy to see that a multichain $C: P_0 \le P_1 \le \cdots \le P_k$ in $\mathcal{P}_n$ has (only) small intervals if $P_i$ is obtained by turning some valleys of $P_{i-1}$ into peaks, for every $i \in [k]$. Saturated $P - \1_n$ chains are those chains with small intervals that have maximum length. In a recent work \cite{MSTT2019}, the authors evaluated the number $f(P)$ of $P - \1_n$ chains with small intervals and minimum length, for a given arbitrary $P \in \mathcal{P}_n$. In this work,
  we study the same subject from another point of view, giving bijective proofs for the main results. In the sequel, we summarize the main definitions and results of \cite{MSTT2019}.

We recall that for every path $P \in \mathcal{P}_n \setminus \{\1_n\}$ the join $\widetilde{P}$ of all elements covering $P$ is called the filling of $P$ and $\widetilde{\1_n} = \1_n$. Obviously, the filling of $P$ is obtained by turning every valley of $P$ into a peak. Moreover, for every path $P \in \mathcal{P}_n$ we define by induction a finite sequence of paths $P^{(i)}$ in $\mathcal{P}_n$ such that $P^{(0)} = P$ and $P^{(i)} = \widetilde{P^{(i-1)}}$ whenever $P^{(i-1)} \ne \mathbf{1}_n$. The number $\delta(P)$ for which $P^{(\delta(P))} = \mathbf{1}_n$ is called {\em degree} of $P$. For every path $P \in \mathcal{P}$ with $P \ne \1_{|P|}$, the degree $\delta(P)$ is given by the formula 
\[ \delta(P) = |P| - 1 - \lv(P).\]
The length of each $P - \1_n$ chain with small intervals and minimum length is equal to $\delta(P)$.
In order to count the number $f(P)$ of these chains, a new kind of multichains of Dyck paths, based on the heights of the valleys of the paths, was introduced: A multichain of Dyck paths $C : \sigma_0 \le \sigma_1 \le \cdots \le \sigma_h$, where $h = \hv(\sigma_0)$, is {\em of type $\V$} iff for every $j \in [h]$ the paths $\sigma_j$, $\sigma_{j-1}$ have the same valleys at every height $\le h - j$. For $a, b \in \mathcal{D}$ with $a \le b$, we denote by $\V(a,b)$ the number of all $a-b$ multichains of type $\V$. Clearly, $\V(a,b) \ne 0$ iff $a, b$ have exactly the same low valleys. 

Since $f(uP) = f(P)$, it is enough to evaluate $f$ only for Dyck prefixes. This has been done by combining the next three results, where the following notation is used: 
\[ \textrm{$I(a) = |[a, u^{|a|/2} d^{|a|/2}]|$, for every $a \in \mathcal{D}$, and $J(P) = |[P, u^{|P|}]|$, for every $P \in \mathcal{P}$.} \]

\begin{proposition}\label{prop:fprefixsuffix} 
The mapping $f$ satisfies the following properties:
\begin{enumerate}[label=\roman{*})]  
\item $f(P_1 P_2) = f(P_1) f(d P_2)$, where $P_1$ is a Dyck suffix and $P_2$ is a Dyck prefix. 
\item $f(d P) = f(P)$, where $P$ is a Dyck prefix that has at least one return point.
\end{enumerate}
\end{proposition}

\begin{proposition}\label{prop:fprimedyck}
For every Dyck path $a \in \mathcal{D}$, we have that
\begin{equation}\label{eq:fprimedyck} 
f(uad) = \sum\limits_{s \ge a} \V(a,s) I(s).  
\end{equation} 
\end{proposition}

\begin{proposition}\label{prop:fdyckprefix}
For every Dyck prefix $P = a_0 u a_1 \cdots u a_k$, $k \ge 0$, $a_l \in \mathcal{D}$, $l \in [0,k]$, we have that
\begin{equation}\label{eq:fdyckprefix1}
f(d u P) = \sum \prod\limits_{l=0}^k \V(a_l, s_l) J(s_0 V_1), 
\end{equation}
where the sum is taken over all finite sequences $(s_l)$, $l \in [0,k]$, of Dyck paths with $s_l \ge a_l$, and over all finite sequences of Dyck prefixes $(V_i)$, $i \in [k+1]$, with $V_i \ge u s_i V_{i+1}$, $i \in [k]$, and $V_{k+1} = \varepsilon$.
\end{proposition}

The above propositions have been proved using induction and an algebraic approach. In this paper, we follow a completely different method. 
Firstly, in section~\ref{section:tableaux}, we exhibit a natural bijection between multichains and shifted tableaux, which is used to enumerate saturated $P - \1_n$ chains and to give a new combinatorial interpretation for the number $f(P)$, using certain shifted increasing tableaux. Then, in section~\ref{section:mapf}, we give bijective proofs of the above propositions by utilizing decompositions of shifted tableaux, thus also enumerating some new classes of shifted tableaux.

\section{Multichains of paths and shifted tableaux}\label{section:tableaux}

A {\em strict partition} of a positive integer $N$ is a (strictly) decreasing sequence $\lambda = (\lambda_1, \lambda_2, \ldots$, $\lambda_m)$ with $\lambda_m > 0$, such that $\sum_{i=1}^m \lambda_i = N$. The integer $N$ is called {\em size} of $\lambda$. The {\em shifted (Ferrers) diagram} of shape $\lambda$ is the set $\{(i,j): i \in [m], i \leq j \leq \lambda_i+i-1\}$
and it is depicted by an array of cells with $m$ rows, 
where each row is indented by one cell to the right with respect to the previous row, and with $\lambda_i$ cells in row $i$. 
A {\em shifted tableau} $T = (t_{i,j})$ (or more simply $T = (t_{ij})$ if there is no ambiguity) of shape $\lambda$ is a filling of the cells of the shifted  diagram of shape $\lambda$ with positive integers, such that the entries along rows and columns are non-decreasing. We denote by $\max(T)$ the maximum element of $T$. 

A shifted tableau $T$ is called {\em strictly increasing} (or simply {\em increasing}) if its entries along rows and columns are strictly increasing. In particular, $T$ is a {\em standard shifted tableau} if every number of the interval $[N]$ appears exactly once in $T$. Increasing tableaux (either left-aligned, or shifted) have been studied by many authors (e.g., see \cite{CTY2014, Pechenik2014, PSV2018, ThomasYong2009}). 

If $T$ is an increasing shifted tableau, it is easy to check that 
\[ t_{ij} \ge i + j - 1 \textrm{ for every $i, j$} \textrm{ and that } 
\max(T) \ge \max\limits_{i \in [m]}(\lambda_i + 2i - 2).  \] 
Clearly, if $t_{ij} = i + j - 1$, then $t_{i^\prime j^\prime} = i^\prime + j^\prime - 1$ for every $i^\prime \le i$ and $j^\prime \le j$.

The shifted diagrams are closely related with the paths. Every path $P$ can be decomposed as
\begin{equation}\label{k_i}
P = u^{k_1}d u^{k_2-k_1} d u^{k_3-k_2} d \cdots u^{k_m-k_{m-1}} d u^{k_{m+1} - k_m},
\end{equation}
where $m=|P|_d$, $k_i$, $i \in [m]$, is the number of $u$'s before the $i$-th downstep of $P$ and $k_{m+1} = |P|_u$. 
In the sequel, we will write $P=(k_i)_{i \in [m+1]}$ to denote the encoding of $P$ by the sequence of these $k_i$'s. Using this encoding, we map every path $P = (k_i)_{i \in [m+1]} \in \mathcal{P}_n$ to a shifted diagram $F(P)$ of shape $\lambda(P) = (\lambda_1, \lambda_2, \ldots, \lambda_m)$, where 
\[
\lambda_i = n - i - k_i + 1, \qquad i \in [m].
\]  
Clearly, the restriction of this mapping to the set of all paths starting with $d$ is a bijection, since in this case $\lambda_1 = n$. 

\begin{figure}[ht]
\begin{center}
\psset{unit=.9em, linewidth=1pt, radius=2pt, labelsep = 2pt}
\begin{tabular}{ccccc}
\begin{tabular}{c}
\psset{unit=0.3in, linewidth=1pt, radius=2pt, labelsep = 2pt}
\begin{pspicture}(0,-3)(6,0)
\psline[linewidth=2pt,linecolor=blue]{-}(6,0)(6,-1)(5,-1)(5,-2)(4,-2)(3,-2)(3,-3)
\psline[]{-}(0,0)(6,-0)
\psline[]{-}(0,-1)(6,-1)
\psline[]{-}(1,-2)(5,-2)
\psline[]{-}(2,-3)(3,-3)
\psline[]{-}(2,0)(2,-3)
\psline[]{-}(3,0)(3,-3)
\psline[]{-}(1,-2)(1,0)
\psline[]{-}(0,-1)(0,0)
\psline[]{-}(4,0)(4,-2)
\psline[]{-}(5,0)(5,-2)
\psline[]{-}(6,0)(6,-1)
\rput(0.5,-0.5){\tiny $(1,1)$}
\rput(1.5,-0.5){\tiny $(1,2)$}
\rput(2.5,-0.5){\tiny $(1,3)$}
\rput(3.5,-0.5){\tiny $(1,4)$}
\rput(4.5,-0.5){\tiny $(1,5)$}
\rput(5.5,-0.5){\tiny $(1,6)$}
\rput(1.5,-1.5){\tiny $(2,2)$}
\rput(2.5,-1.5){\tiny $(2,3)$}
\rput(3.5,-1.5){\tiny $(2,4)$}
\rput(2.5,-2.5){\tiny $(3,3)$}
\rput(4.5,-1.5){\tiny $(2,5)$}
\rput(3,1){\tiny $(i,j)$}
\end{pspicture}
\end{tabular} 
& $\leftrightarrow$ &
\begin{tabular}{c}
\psset{unit=0.212132in, linewidth=1pt, radius=2pt, labelsep = 2pt}
\begin{pspicture}(0,-1)(7,6)
\psline[linewidth=2pt,linecolor=blue]{-}(0,0)(1,-1)(2,0)(3,-1)(4,0)(5,1)(6,0)
\psline[linestyle=dotted]{-}(0,0)(6,6)
\psline[linestyle=dotted]{-}(0,0)(1,-1)
\psline[linestyle=dotted]{-}(1,-1)(7,5)
\psline[linestyle=dotted]{-}(7,5)(6,6)
\psline[linestyle=dotted]{-}(3,-1)(7,3)
\psline[linestyle=dotted]{-}(6,0)(7,1)
\psline[linestyle=dotted]{-}(3,-1)(1,1)
\psline[linestyle=dotted]{-}(4,0)(2,2)
\psline[linestyle=dotted]{-}(6,0)(3,3)
\psline[linestyle=dotted]{-}(7,1)(4,4)
\psline[linestyle=dotted]{-}(7,3)(5,5)
\psdot[linewidth=0.07em](1,-1)
\psdot[linewidth=0.07em](2,0)
\psdot[linewidth=0.07em](3,-1)
\psdot[linewidth=0.07em](4,0)
\psdot[linewidth=0.07em](5,1)
\psdot[linewidth=0.07em](6,0)
\psdot[linewidth=0.07em](3,1)
\psdot[linewidth=0.07em](4,2)
\psdot[linewidth=0.07em](5,3)
\psdot[linewidth=0.07em](6,4)
\psdot[linewidth=0.07em](6,2)
\rput(6,5){\tiny $(6,4)$}
\rput(5,4){\tiny $(5,3)$}
\rput(4,3){\tiny $(4,2)$}
\rput(3,2){\tiny $(3,1)$}
\rput(6,3){\tiny $(6,2)$}
\rput(2,1){\tiny $(2,0)$}
\rput(1,0){\tiny $(1,-1)$}
\rput(5,2){\tiny $(5,1)$}
\rput(4,1){\tiny $(4,0)$}
\rput(6,1){\tiny $(6,1)$}
\rput(3,0){\tiny $(3,-1)$}
\rput(3,6){\tiny $(n+i-j,n-i-j)$}
\end{pspicture} 
\end{tabular}\\
\end{tabular} 
\end{center}
\caption{The mapping of the cells $(i,j)$ of $F(P)$ to the lattice points $(n+i-j, n-i-j)$ weakly above $P = d u d u^2 d$} \label{fig:ferrers}
\end{figure}
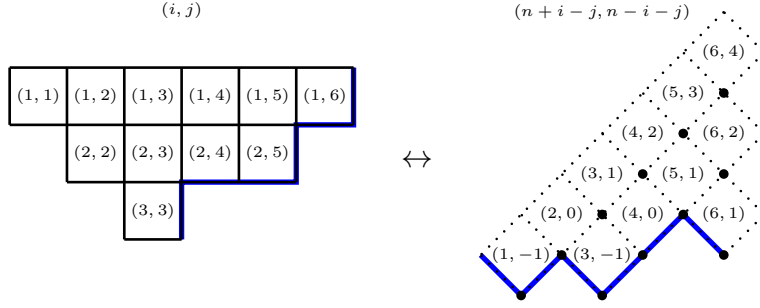

We can easily check that the mapping $(i,j) \mapsto (n+i-j, n-i-j)$ maps bijectively the cells of the shifted diagram $F(P)$ to the lattice points that are weakly above $P$, strictly below $\1_n$ and have coordinates of the same parity. In particular, the cells $(i,j)$, $j \ge m$, that have their south-eastern corner on the border of the shifted diagram correspond to the lattice points of the path $P$ (excluding the initial point $(0,0)$), e.g., see Figure \ref{fig:ferrers} for the path $P = d u d u^2 d$, colored blue. Thus, the path $P$ is obtained by rotation and reflection of the staircase path starting from the north-eastern corner of the cell $(1,n)$ and passing through these corner points. 

Next, some easily verified facts concerning the shifted diagrams are stated.

\begin{remark}\label{rem:mainremark} \ 

\begin{enumerate}[label=\roman{*})]

\item\label{rem:r1} For every $P, Q \in \mathcal{P}_n$, we have that $P \le Q$ iff $F(Q)$ is a subdiagram of $F(P)$.

\item\label{rem:r2} Let $P \le Q$ and a cell $(i,j)$ of $F(Q)$; if the point $A(n+i-j, n-i-j)$ belongs to the path $P$, then it also belongs to the path $Q$. 
\item\label{rem:r3} A cell $(i,j)$ of $F(P)$ corresponds to a low valley point of $P$ iff $i + j - 1 = \delta(P)$. 

\item\label{rem:r4} If the corresponding points of two cells $(i,j)$, $(i^\prime, j^\prime)$ of $F(P)$ are path-connected (i.e., there exists a lattice path connecting them) and $j^\prime - i^\prime < j - i$, then $j^\prime \le j$. 
\end{enumerate}

\end{remark}

In the sequel, we define a mapping between $P - \1_n$ multichains $C$ of length $k$ and shifted tableaux $T$ of shape $\lambda(P)$ with $\max(T) \le k$. Given a multichain $C : P_0 = P \le P_1 \le \cdots \le P_k = \1_n$ and a cell $(i,j) \in F(P)$, there exists a unique $\xi \in [0,k-1]$ such that the point $A(n+i-j,n-i-j)$ lies weakly above $P_\xi$ and strictly below $P_{\xi+1}$. Then, we define $T = (t_{ij})$ with $t_{ij} = k - \xi$, i.e., 
\begin{equation}\label{eq:bijectiontheta} 
t_{ij} = k - \xi \Leftrightarrow \textrm{ $h_{n+i-j}(P_\xi) \le n - i - j < h_{n+i-j}(P_{\xi+1})$}.
\end{equation}
Conversely, in order to construct a unique multichain from a shifted tableau, we must know the length $k$ of the desired multichain. More precisely, given a pair $(T,k)$ where $T = (t_{ij})$ is a shifted tableau of shape $\lambda(P)$ with $\max(T) \le k$, we recover the original chain $C = (P_\xi)$, $0 \le \xi \le k$, where $P_\xi$ is defined by the staircase path that is determined by the boundary between the cells with numbers less than or equal to $k-\xi$ and the cells with numbers greater than $k-\xi$.

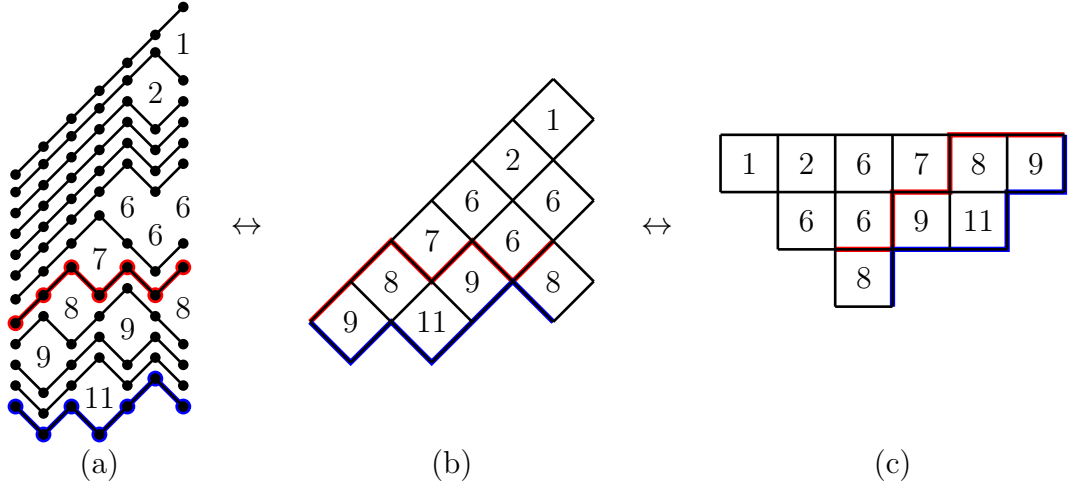
\begin{figure}[ht]
\begin{center}
\psset{unit=.9em, linewidth=1pt, radius=2pt, labelsep = 2pt}
\begin{tabular}{ccccc}
\begin{tabular}{c}
\begin{pspicture}(0,0)(6,6)
\psline[showpoints=true]{-}(0,0)(1,1)(2,2)(3,3)(4,4)(5,5)(6,6)
\end{pspicture} \\[-4.3em]
\begin{pspicture}(0,0)(6,5)
\psline[showpoints=true]{-}(0,0)(1,1)(2,2)(3,3)(4,4)(5,5)(6,4)
\rput(6,5.3){$1$}
\end{pspicture} \\[-3.3em]
\begin{pspicture}(0,0)(6,4)
\psline[showpoints=true]{-}(0,0)(1,1)(2,2)(3,3)(4,4)(5,3)(6,4)
\rput(5,4.3){$2$}
\end{pspicture} \\[-3.3em]
\begin{pspicture}(0,0)(6,4)
\psline[showpoints=true]{-}(0,0)(1,1)(2,2)(3,3)(4,4)(5,3)(6,4)
\end{pspicture} \\[-3.3em]
\begin{pspicture}(0,0)(6,4)
\psline[showpoints=true]{-}(0,0)(1,1)(2,2)(3,3)(4,4)(5,3)(6,4)
\end{pspicture} \\[-3.3em]
\begin{pspicture}(0,0)(6,4)
\psline[showpoints=true]{-}(0,0)(1,1)(2,2)(3,3)(4,4)(5,3)(6,4)
\end{pspicture} \\[-2.3em]
\begin{pspicture}(0,0)(6,3)
\psline[showpoints=true]{-}(0,0)(1,1)(2,2)(3,3)(4,2)(5,1)(6,2)
\rput(5,2.3){$6$}
\rput(4,3.3){$6$}
\rput(6,3.3){$6$}
\end{pspicture} \\[-2.3em]
\begin{pspicture}(0,0)(6,3)
\psline[showpoints=true,linecolor=red,linewidth=2pt]{-}(0,0)(1,1)(2,2)(3,1)(4,2)(5,1)(6,2)
\psline[showpoints=true]{-}(0,0)(1,1)(2,2)(3,1)(4,2)(5,1)(6,2)
\rput(3,2.3){$7$}
\end{pspicture} \\[-1.5em]
\begin{pspicture}(0,0)(6,2)
\psline[showpoints=true]{-}(0,0)(1,1)(2,0)(3,1)(4,2)(5,1)(6,0)
\rput(2,1.3){$8$}
\rput(6,1.3){$8$}
\end{pspicture} \\[-0.6em]
\begin{pspicture}(0,-1)(6,1)
\psline[showpoints=true]{-}(0,0)(1,-1)(2,0)(3,1)(4,0)(5,1)(6,0)
\rput(1,0.3){$9$}
\rput(4,1.3){$9$}
\end{pspicture} \\[-1.5em]
\begin{pspicture}(0,-1)(6,1)
\psline[showpoints=true]{-}(0,0)(1,-1)(2,0)(3,1)(4,0)(5,1)(6,0)
\end{pspicture} \\[-1.5em]
\begin{pspicture}(0,-1)(6,1)
\psline[showpoints=true,linecolor=blue,linewidth=2pt]{-}(0,0)(1,-1)(2,0)(3,-1)(4,0)(5,1)(6,0)
\psline[showpoints=true]{-}(0,0)(1,-1)(2,0)(3,-1)(4,0)(5,1)(6,0)
\rput(3,0.3){$11$}
\end{pspicture} 
\end{tabular}
& $\leftrightarrow$ &
\begin{tabular}{c}
\psset{unit=0.212132in, linewidth=1pt, radius=2pt, labelsep = 2pt}
\begin{pspicture}(0,-1)(7,6)
\psline[linewidth=2pt,linecolor=red]{-}(0,0)(2,2)(3,1)(4,2)(5,1)(6,2)
\psline[linewidth=2pt,linecolor=blue]{-}(0,0)(1,-1)(2,0)(3,-1)(4,0)(5,1)(6,0)
\psline[]{-}(0,0)(6,6)
\psline[]{-}(0,0)(1,-1)
\psline[]{-}(1,-1)(7,5)
\psline[]{-}(7,5)(6,6)
\psline[]{-}(3,-1)(7,3)
\psline[]{-}(6,0)(7,1)
\psline[]{-}(3,-1)(1,1)
\psline[]{-}(4,0)(2,2)
\psline[]{-}(6,0)(3,3)
\psline[]{-}(7,1)(4,4)
\psline[]{-}(7,3)(5,5)
\rput(6,5){$1$}
\rput(5,4){$2$}
\rput(4,3){$6$}
\rput(3,2){$7$}
\rput(6,3){$6$}
\rput(2,1){$8$}
\rput(1,0){$9$}
\rput(5,2){$6$}
\rput(4,1){$9$}
\rput(6,1){$8$}
\rput(3,0){$11$}
\end{pspicture} 
\end{tabular}
& $\leftrightarrow$ &
\begin{tabular}{c}
\psset{unit=0.3in, linewidth=1pt, radius=2pt, labelsep = 2pt}
\begin{pspicture}(0,-3)(6,0)
\psline[linewidth=2pt,linecolor=red]{-}(6,0)(4,0)(4,-1)(3,-1)(3,-2)(2,-2)
\psline[linewidth=2pt,linecolor=blue]{-}(6,0)(6,-1)(5,-1)(5,-2)(4,-2)(3,-2)(3,-3)
\psline[]{-}(0,0)(6,-0)
\psline[]{-}(0,-1)(6,-1)
\psline[]{-}(1,-2)(5,-2)
\psline[]{-}(2,-3)(3,-3)
\psline[]{-}(2,0)(2,-3)
\psline[]{-}(3,0)(3,-3)
\psline[]{-}(1,-2)(1,0)
\psline[]{-}(0,-1)(0,0)
\psline[]{-}(4,0)(4,-2)
\psline[]{-}(5,0)(5,-2)
\psline[]{-}(6,0)(6,-1)
\rput(0.5,-0.5){$1$}
\rput(1.5,-0.5){$2$}
\rput(2.5,-0.5){$6$}
\rput(3.5,-0.5){$7$}
\rput(1.5,-1.5){$6$}
\rput(4.5,-0.5){$8$}
\rput(5.5,-0.5){$9$}
\rput(2.5,-1.5){$6$}
\rput(3.5,-1.5){$9$}
\rput(2.5,-2.5){$8$}
\rput(4.5,-1.5){$11$}
\end{pspicture}
\end{tabular} \\
(a) & & (b) & & (c)
\end{tabular} 
\end{center}
\caption{A $d u d u^2 d - \1_6$ multichain and the corresponding shifted tableau} \label{fig:multichainstoshiftedtableaux}
\end{figure}

The above procedure determines a bijection between multichains of paths and shifted tableaux, that can be described practically as in Figure \ref{fig:multichainstoshiftedtableaux}, where the integers in Figure \ref{fig:multichainstoshiftedtableaux}(a) count the number of paths in the chain that lie above a given lattice point, so that they generate the tableau of Figure \ref{fig:multichainstoshiftedtableaux}(b), which by rotation and reflection gives the corresponding shifted tableau of Figure \ref{fig:multichainstoshiftedtableaux}(c). 

In the following we give some characterizations for special kinds of multichains.

Firstly, for $t \in [k]$ we have that $P_{k-t} \ne P_{k-t+1}$ iff there exists a point $(n+i-j,n-i-j)$ which lies weakly above $P_{k-t}$ and strictly below $P_{k-t+1}$, i.e., iff there exists a cell $(i,j)$ of $T$ such that $t_{ij} = t$. In particular, $C$ is a chain iff the entries of $T$ form the interval $[k]$.

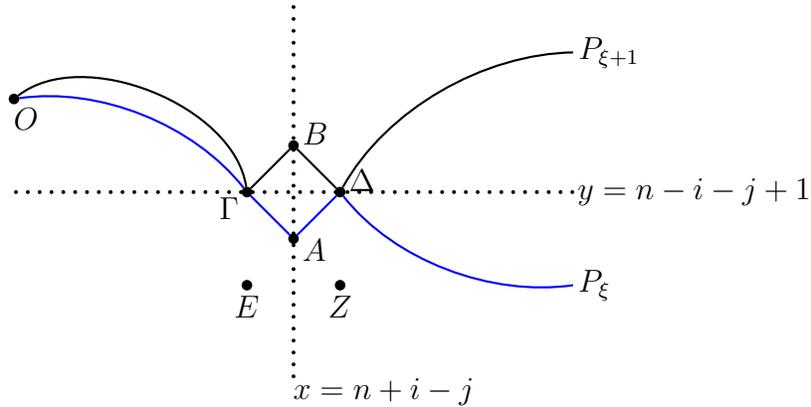
\begin{figure}[ht]
\begin{center}
\psset{unit=1.5em, radius=2pt,labelsep=2pt}
\begin{pspicture}(1,-4)(18,4)

\Cnode*(6,0){C} \nput[]{220}{C}{$\Gamma$}
\Cnode*(8,0){D} \nput[]{30}{D}{$\Delta$}
\Cnode*(7,1){B} \nput[]{30}{B}{$B$}
\Cnode*(7,-1){A} \nput[]{-30}{A}{$A$}
\Cnode*(6,-2){E} \nput[]{-90}{E}{$E$}
\Cnode*(8,-2){Z} \nput[]{-90}{Z}{$Z$}


\psline[linestyle=dotted, linewidth=1.5pt](7,-4)(7,4) \rput[lt](7,-4){$x=n+i-j$}
\psline[linestyle=dotted, linewidth=1.5pt](1,0)(13,0) \rput[l](13.1,0){$y=n-i-j+1$}

\Cnode*(1,2){O} \nput[]{-60}{O}{$O$}
\pnode(13,3){P1} \nput[]{0}{P1}{$P_{\xi+1}$}
\ncarc[arcangle=60]{O}{C} \ncarc[arcangle=30]{D}{P1}
\ncline[]{C}{B} \ncline[]{B}{D}


\pnode(13,-2){P2} \nput[]{0}{P2}{$P_{\xi}$}
\ncarc[arcangle=30, linecolor=blue]{O}{C} \ncarc[arcangle=-30, linecolor=blue]{D}{P2}
\ncline[linecolor=blue]{C}{A} \ncline[linecolor=blue]{A}{D}

\end{pspicture}
\end{center}
\caption{Two consecutive paths of the multichain generated by an increasing tableau}\label{fig:rombus}
\end{figure}

Secondly, if $C$ has small intervals and $t_{ij} = k - \xi$, by relation \eqref{eq:bijectiontheta} we have that 
$h_{n+i-j}(P_{\xi+1}) = h_{n+i-j}(P_{\xi}) + 2 = n - i -j + 2$, and hence the point $A(n+i-j,n-i-j)$ is a valley of $P_\xi$ and the point $B(n+i-j,n-i-j+2)$ is a peak of $P_{\xi+1}$. It follows that the points $\Gamma(n+(i-1)-j,n-(i-1)-j)$ and $\Delta(n+i-(j-1), n-i-(j-1))$ lie on the path $P_{\xi+1}$. Hence, $t_{i-1,j}, t_{i,j-1} < k - \xi = t_{ij}$, i.e., $T$ is increasing. Conversely, assume that $T$ is increasing and let $\xi \in [0,k-1]$ and $x \in [n]$, such that $h_x(P_\xi) < h_x(P_{\xi+1})$. Let $(i,j) \in F(P)$ such that $n+i-j = x$ and $n-i-j = h_x(P_{\xi})$; it follows that $t_{ij} = k-\xi$. Then, since $T$ is increasing we have that $t_{i-1,j}, t_{i,j-1} < k - \xi = t_{ij}$, so that the points $\Gamma(n+(i-1)-j,n-(i-1)-j)$ and $\Delta(n+i-(j-1), n-i-(j-1))$ lie on $P_{\xi+1}$, and so the point $B(n+i-j,n-i-j+2)$ is a peak of $P_{\xi+1}$ with $h_x(P_{\xi+1}) = h_x(P_{\xi}) + 2$. Moreover, since $t_{i,j+1}, t_{i+1,j} > k - \xi$, the path $P_\xi$ passes strictly above the points $E(n+i-(j+1),n-i-(j+1))$ and $Z(n+(i+1)-j,n-(i+1)-j)$; hence, $P_\xi$ passes from the points $\Gamma$ and $\Delta$, so that the point $A$ is a valley of $P_\xi$. This shows that $C$ has small intervals (see Figure \ref{fig:rombus}).

Finally, since the length of every saturated $P - \1_n$ chain is equal to  $l(P,\1_n) = \lambda_1 + \lambda_2 + \cdots + \lambda_m$, using the first characterization, we obtain that $C$ is a saturated chain iff the entries of $T$ form the interval $[\lambda_1 + \lambda_2 + \cdots + \lambda_n]$, i.e., iff $T$ is standard. 

Summarizing, we obtain the following result.

\begin{proposition}\label{prop:bijectiontheta}
If $P \in \mathcal{P}_n$ starts with $d$ and $k \in \mathbb{N}^*$ then, there exists a bijection between the $P-\1_n$ multichains $C$ of length $k$ and the shifted tableaux $T$ of shape $\lambda(P)$ with $\max(T) \le k$. In particular,
\begin{enumerate}[label=\roman{*})]
\item $C$ is a chain iff the entries of $T$ form the interval $[k]$.

\item $C$ has small intervals iff $T$ is increasing.
 
\item $C$ is a saturated chain iff $T$ is standard.
\end{enumerate}
\end{proposition}

We note that from the construction of the bijection of the previous proposition, it follows that the number $k - \max(T) + 1$ counts the members of the chain which are equal to $P$. Furthermore, if a path $P \in \mathcal{P}_n$ does not necessarily start with $d$ and the length of its first ascent is equal to $\nu$, then by applying Proposition \ref{prop:bijectiontheta} for the path obtained by deleting the first ascent of $P$, the tableau that we obtain has a first row of length $n - \nu$, giving the following result.

\begin{corollary}
The number of multichains $C$ of the form $P_0 = P_1 = \cdots = P_{\mu-1} < P_{\mu} \le \cdots \le P_k = \1_n$, where $P_0 \in \mathcal{P}_n$ and it has length of the first ascent equal to $\nu$, is equal to the number of shifted tableaux $T$ with length of the first row equal to $n - \nu$ and $\max(T) = k - \mu + 1$. 
\end{corollary}

Clearly, since the minimum length of every $P - \1_n$ chain with small intervals is equal to $\delta(P)$, by Proposition~\ref{prop:bijectiontheta} it follows that \[ \begin{aligned} f(P) =  \# \textrm{ increasing shifted tableaux $T$ of shape $\lambda(P)$} \\ \textrm{such that the entries of $T$ form the interval $[\delta(P)]$}. \end{aligned} \]
Note that in this case, the condition ``the entries of $T$ form the interval $[\delta(P)]$'' can be replaced by the weaker condition ``$\max(T) = \delta(P)$''. Indeed, if $(i,j)$ is a cell of $F(P)$ which corresponds to a low valley point of $P$, then $t_{ij} = i + j - 1 = \delta(P)$ (see Remark \ref{rem:mainremark}\ref{rem:r3}), and since $T$ is increasing we deduce that there exist $\delta(P)$ different entries of $T$ with values less than or equal to $\delta(P)$, giving that the entries of $T$ form the interval $[\delta(P)]$.

Summarizing, we obtain the following result.

\begin{corollary}\label{cor:basic}
For every path $P \in \mathcal{P}_n$ that starts with $d$, we have that
\[ f(P) = \# \textrm{ increasing shifted tableaux $T$ of shape $\lambda(P)$, with $\max(T) = \delta(P)$}. \]
\end{corollary}

Finally, given any path $P = (k_i)_{i \in [m]}$ in $\mathcal{P}_n \setminus \{\1_n\}$, by applying Proposition~\ref{prop:bijectiontheta} iii) for the path obtained by deleting the first ascent of $P$, we deduce that the number of saturated $P - \1_n$ chains is equal to the number of standard shifted tableaux of shape $\lambda = (\lambda_1, \lambda_2, \ldots, \lambda_m)$, where $\lambda_i = n - i - k_i + 1$, $i \in [m]$. Then, using the well-known formula for standard shifted tableaux (e.g., see \cite[p. 267]{MacDonald}) we obtain the following result.

\begin{corollary}
For every path $P = (k_i)_{i \in [m]}$ in $\mathcal{P}_n \setminus \{ \1_n\}$, 
the number of saturated $P - \1_n$ chains is equal to
\begin{equation*}
\dfrac{(\lambda_1 + \lambda_2 + \cdots + \lambda_m)!}{\lambda_1 ! \lambda_2 ! \cdots \lambda_m!} \prod\limits_{i < j} \dfrac{\lambda_i - \lambda_j}{\lambda_i + \lambda_j},
\end{equation*}
where $\lambda_i = n - i - k_i + 1$, $i \in [m]$.
\end{corollary}

\section{Bijective proofs}\label{section:mapf}

In this section we prove Propositions 1, 2 and 3 via increasing shifted tableaux.

\begin{proof}[\bf Proof of Proposition \ref{prop:fprefixsuffix}] \

\noindent{i)} Without loss of generality we may assume that $P_1$ starts with $d$, $P_1 \ne d$ and $P_2 \ne \varepsilon$. Let $n_1 = |P_1|$, $n_2 = |P_2|$, $m_1 = |P_1|_d$, $m_2 = |P_2|_d$ and $\lambda(P_1) = (\lambda_i^{1})_{i \in [m_1]}$, $\lambda(P_2) = (\lambda_i^{2})_{i \in [m_2]}$. Clearly, $\lambda(d P_2) = (\lambda^{2}_1 + 1, \lambda^{2}_1, \lambda^{2}_2, \ldots, \lambda^{2}_{m_2})$ and $\lambda(P_1 P_2) = (\lambda_1^{1} + n_2, \lambda_2^{1} + n_2, \ldots, \lambda_{m_1}^{1} + n_2, \lambda_1^{2}, \lambda_2^{2}, \ldots, \lambda_{m_2}^{2})$.

In view of Corollary \ref{cor:basic}, it is enough to construct a bijection between increasing shifted tableaux $T$ of shape $\lambda(P_1 P_2)$ with $\max(T) = \delta(P_1 P_2) = 2m_1 + n_2 - 1$, and pairs ($T_1$, $T_2$) of increasing shifted tableaux of shapes $\lambda(P_1)$ and $\lambda(d P_2)$ with maximum elements $\max(T_1) = \delta(P_1) = 2m_1 - 1$ and $\max(T_2) = \delta(d P_2) = n_2 + 1$ respectively.

Indeed, for $T = (t_{ij})$ we define 
$T_1 = (t_{ij}^1)$ and $T_2 = (t_{ij}^2)$ 
as follows:
\[ t^1_{ij} = t_{i, j + n_2} - n_2, \textrm{ where } (i,j) \in F(P_1) \]
and
\[ t^2_{ij} = t_{i+m_1-1, j+m_1-1} - 2(m_1-1), \textrm{ where } (i,j) \in F(dP_2).\]

Clearly, since $T$ is an increasing tableau, $t_{m_1, m_1 + n_2} = 2m_1 + n_2 - 1$. Hence, $t_{ij} = i + j - 1$ for every $i \le m_1$ and $j \le m_1 + n_2$. It follows that $t^1_{11} = 1 = t^2_{11}$. Since the entries of $T_1$, $T_2$ are increasing along rows and columns, we obtain that $T_1$, $T_2$ are increasing shifted tableaux of shapes $\lambda(P_1)$ and $\lambda(d P_2)$ respectively, with
\[ \max(T_1) = t^1_{m_1 m_1} = t_{m_1, m_1 + n_2} - n_2 = 2m_1 -1, \]
and
\[ \max(T_2) = t^2_{1, n_2 +1} = t_{m_1, m_1 + n_2} - 2(m_1 - 1) = n_2 + 1.\]

From the above construction we see that $T$ can be decomposed into three parts. Two of them generate $T_1$, $T_2$ by translation, and for the third (north-western) part the corresponding cells $(i,j)$ have $t_{ij} = i + j - 1$ (see Figure \ref{fig:decompi}).

\begin{figure}[ht]
\begin{center}
\includegraphics{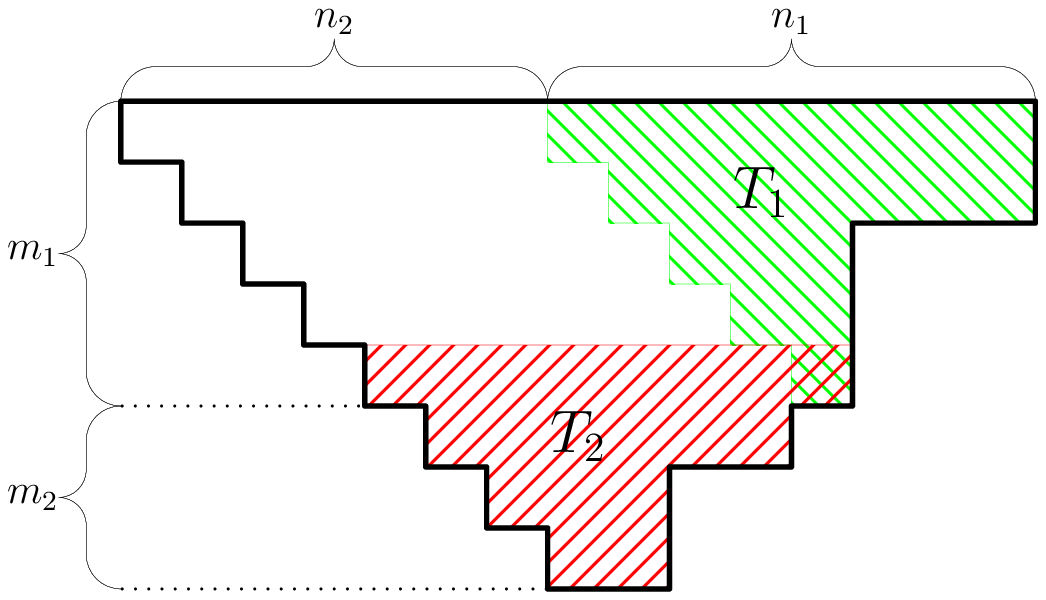}
\end{center}
\caption{A decomposition of a shifted tableau corresponding to a product of paths}\label{fig:decompi}
\end{figure}

Clearly, this procedure is reversible, so that we obtain the desired bijection. 

\noindent{ii)} Let $P = u^k Q$, where $k > 0$ and $Q$ starts with $d$. It is enough to prove that $f(d P) = f(Q)$.  In view of Corollary \ref{cor:basic}, it is enough to construct a bijection between increasing shifted tableaux $T = (t_{ij})$ of shape $\lambda( dP ) = (\lambda_1, \ldots, \lambda_{m+1})$, where $m = |P|_d$ with $\max(T) = \delta( d P) = |P| + 1$, and increasing shifted tableaux $T^\prime = (t^\prime_{ij})$ of shape $\lambda(Q) = (\lambda_2, \ldots, \lambda_{m+1})$ with $\max(T^\prime) = \delta(Q) = |P| -1$.

A return point of $P$ corresponds to a cell $(r,c)$ of $T$ such that $t_{rc} = r + c - 1 = |P| + 1$ and $2 \le r \le c$, (see Remark \ref{rem:mainremark}\ref{rem:r3}). Thus, $t_{ij} = i + j - 1$ whenever $i \le r$ and $j \le c$ and, in particular, $t_{2,2} = 3$. Moreover, since the first row of $T$ has $|P| + 1 = \max(T)$ cells, it follows that $t_{1j} = j$, for all $j \in [|P| + 1]$. Then, $T^\prime$ is obtained by deleting the first row of $T$ and by subtracting 2 from each of the remaining elements, i.e., 
\[ t^\prime_{ij} = t_{i +1, j + 1} - 2. \]
Clearly, $T^\prime$ satisfies the required conditions. On the other hand, by adding 2 to each element of $T^\prime$, and by adding an extra first row of length $|P|+1$ consisting of the interval $[|P|+1]$, we recover $T$ (see Figure \ref{fig:decompii}). 
\begin{figure}[ht]
\begin{center}
\includegraphics{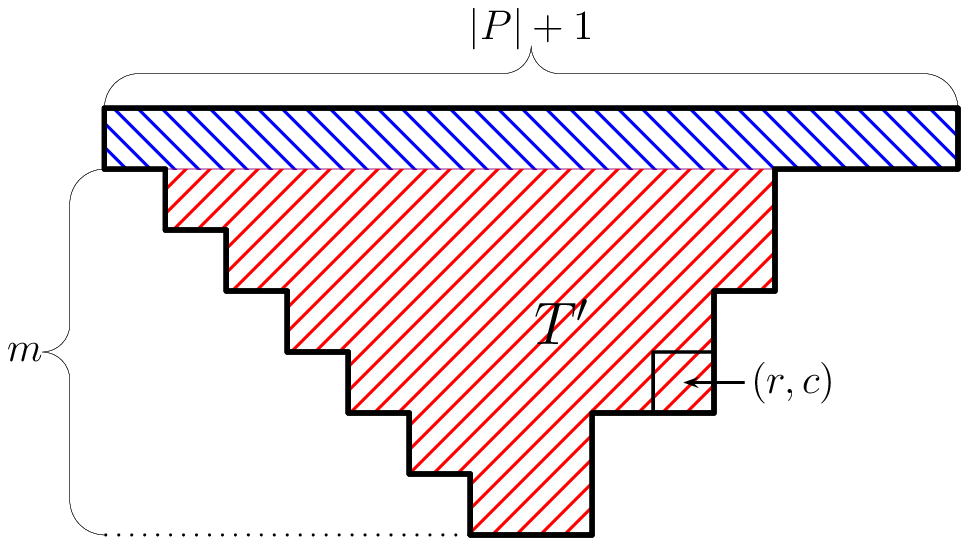}
\end{center}
\caption{A decomposition of a shifted tableau corresponding to a path starting with $d$}\label{fig:decompii}
\end{figure}
\end{proof}

\begin{proof}[\bf Proof of Proposition \ref{prop:fprimedyck}] Set $a = (k_i)_{i \in [m+1]} \in \mathcal{D}_m$ and let $P$ be the suffix of $a$ starting at the first $d$ of $a$, i.e., $a = u^{k_1} P$. Then, $P d$ is written in the form $P d = (k_i - k_1)_{i \in [m+2]}$, where $k_{m+2} = k_{m+1} = k_m = m$.

Applying Corollary \ref{cor:basic} for the path $P d$, we obtain that
\begin{equation}\label{eq:eqb1}
\begin{aligned}
f(uad) = f(Pd) = & \# \textrm{ increasing shifted tableaux $T$ of shape} \\ 
& \textrm{$\lambda(Pd) = (\lambda_1, \lambda_2, \ldots, \lambda_{m+1}$), with $\max(T) = \delta(Pd) = 2m+1$,}
\end{aligned}
\end{equation}
where $\lambda_i = 2m + 2 - i - k_i$, $i \in [m+1]$.

We first show that for any such tableau $T = (t_{ij})$ we have that
\begin{equation}\label{eq:tijfixedpart}
t_{ij} = i + j - 1, \textrm{ for $j \le m + 1$,} 
\end{equation}
and
\begin{equation}\label{eq:tijupperbound}
t_{ij} \le i + j + \hv(a), \textrm{ for $(i,j) \in F(Pd)$}. 
\end{equation}
Indeed, since $t_{m+1,m+1} = 2m + 1$, formula \eqref{eq:tijfixedpart} is a direct consequence of the increasing property of $T$. Moreover, we note that if \eqref{eq:tijupperbound} holds for either one of the cells $(i+1,j)$, $(i,j+1)$ then, using the increasing property of $T$, it must also hold for the cell $(i,j)$. Thus, it is enough to show \eqref{eq:tijupperbound} for every cell $(i,j)$ which is the rightmost cell in its row as well as the lowest cell in its column, i.e., $j = i + \lambda_i - 1 = 2m + 1 - k_i$ and $k_{i+1} > k_i$. Clearly, since each one of these cells corresponds to a valley point of $a$, we have that $k_i - i \le \hv(a)$, and since $\max(T) = 2m + 1$, we conclude that
\[ t_{i,2m+1-k_i} \le i + (2m + 1 - k_i) + \hv(a), \]
so that \eqref{eq:tijupperbound} holds.

On the other hand, we note that the sum of the right-hand side of formula \eqref{eq:fprimedyck} counts the multichains $(\sigma_i)_{i \in [0,k]}$, where $\sigma_ 0 = a$, $\sigma_{k-1} = u^m d^m$, $\sigma_k = u^{2m}$, such that the sub-multichain $\sigma_0 \le \sigma_1 \le \cdots \le \sigma_{k-3}$ is of type $\V$, and $k = \hv(a) + 3$, or equivalently it counts the multichains $(s_i)_{i \in [0,k]}$, where each $s_i$ is obtained from $\sigma_i$ by deleting the first $k_1$ upsteps of $\sigma_i$. 

Using the bijection of Proposition \ref{prop:bijectiontheta} for the path $P$, we deduce that each one of these multichains $(s_i)_{i \in [0,k]}$ corresponds bijectively to a shifted tableau $V = (v_{ij})$ of shape $\lambda(P) = (\lambda_1-1, \lambda_2-1, \ldots, \lambda_m-1)$, with $\max (V) \le k$. Furthermore, using relation \eqref{eq:bijectiontheta} we have that
\begin{align}\label{eq:vij} 
v_{ij} = k - \xi & \Leftrightarrow h_{2m - k_1 + i - j} (s_\xi) \le 2m - k_1 - i - j < h_{2m - k_1 + i - j} (s_{\xi+1}) \nonumber \\
& \Leftrightarrow h_{2m + i - j} (\sigma_\xi) \le 2m - i - j < h_{2m + i - j}(\sigma_{\xi+1}).
\end{align}


We will show that
\begin{equation}\label{eq:vijlowerbound}
v_{ij} = 1 \textrm{ iff } j \le m,
\end{equation}
and
\begin{equation}\label{eq:vijupperbound}
v_{ij} \le 2m + 3 - i - j, \textrm{ for $(i,j) \in F(P)$.}
\end{equation}
Firstly, since $\sigma_{k-1} = u^m d^m$, we have that
\[ h_{2m + i - j}(\sigma_{k-1}) = 
\begin{cases} j - i, & \textrm{ if $j \le i + m$;} \\ 
2m + i - j, & \textrm{ if $j > i + m$;} 
\end{cases}\]
so, by applying relation \eqref{eq:vij} for $\xi = k-1$ we obtain relation \eqref{eq:vijlowerbound}.

Clearly, by the non-decreasing property of the tableau $V$, it is enough to show relation \eqref{eq:vijupperbound} for $j = i + \lambda_i - 2$ and $i + 1 + \lambda_{i+1} - 2 < i + \lambda_i - 2$, i.e., for $j = 2m - k_i$ and $k_i < k_{i+1}$. We note that, in this case the cell $(i,j)$ corresponds to the point $(k_i+i, k_i-i)$ which is a valley point of $a$, so that $k_i - i \in [0,k-3]$. Now, since the multichain $\sigma_0 \le \sigma_1 \le \cdots \le \sigma_{k-3}$ is of type $\V$, the paths $a$ and $\sigma_{k-3 - (k_i - i)}$ coincide up to height $k_i - i$, so that
\[ h_{k_i + i} (\sigma_{k-3-(k_i-i)}) = h_{k_i + i}(a) = k_i - i. \]
Then, if $v_{i,2m-k_i} = k - \xi$, by relation \eqref{eq:vij} we have that
$k_i - i < h_{k_i + i}(\sigma_{\xi+1})$, so that by the above equality we obtain that $k - 3 - (k_i-i) < \xi + 1$. Hence, 
\[ v_{i,2m-k_i} = k - \xi \le 3 + k_i - i = 2m + 3 - i - (2m-k_i).\]  

Conversely, assuming that the tableau $V$ satisfies relations \eqref{eq:vijlowerbound} and \eqref{eq:vijupperbound}, we can easily check that the multichain $\sigma_0 \le \sigma_1 \le \cdots \le \sigma_{k-3}$ is of type $\V$, and $\sigma_{k-1} = u^m d^m$, which gives that
\begin{equation}\label{eq:eqb2}
\begin{aligned}
\sum\limits_{s \ge a} v(a,s) I(s) = \# \textrm{ shifted tableaux $V$ of shape $\lambda = (\lambda_1-1, \lambda_2-1, \ldots, \lambda_{m}-1$)}  \\ 
\textrm{ such that the entries $v_{ij}$ of $V$ satisfy relations \eqref{eq:vijlowerbound} and \eqref{eq:vijupperbound}.}
\end{aligned}
\end{equation}

In view of formulas \eqref{eq:eqb1} and \eqref{eq:eqb2}, for the justification of formula \eqref{eq:fprimedyck} it is enough to give a bijection between the tableaux $T = (t_{ij})$ and $V = (v_{ij})$.

For this, set
\[ v_{ij} = \begin{cases} 
1, & \textrm{ if $j \le m$;} \\
t_{i,j+1} + 2 - i - j, & \textrm{ if $j \ge m + 1$.} 
\end{cases} \]
It is easy to check that $V$ is a shifted tableau, with $v_{ij} \ge 2$ iff $j \ge m+1$, so that relation \eqref{eq:vijlowerbound} holds. Since $\max(T) = 2m+1$, we obtain that relation \eqref{eq:vijupperbound} holds too, and by relation \eqref{eq:tijupperbound}, we obtain automatically that $\max(V) \le k$. For the reverse, set
\[ t_{ij} = \begin{cases}
i + j - 1 & \textrm{ if $j \le m+1$;} \\
v_{i,j-1} + i + j -3 & \textrm{ if $j \ge m+2$,}
\end{cases} \]
giving the required bijection.
\end{proof}

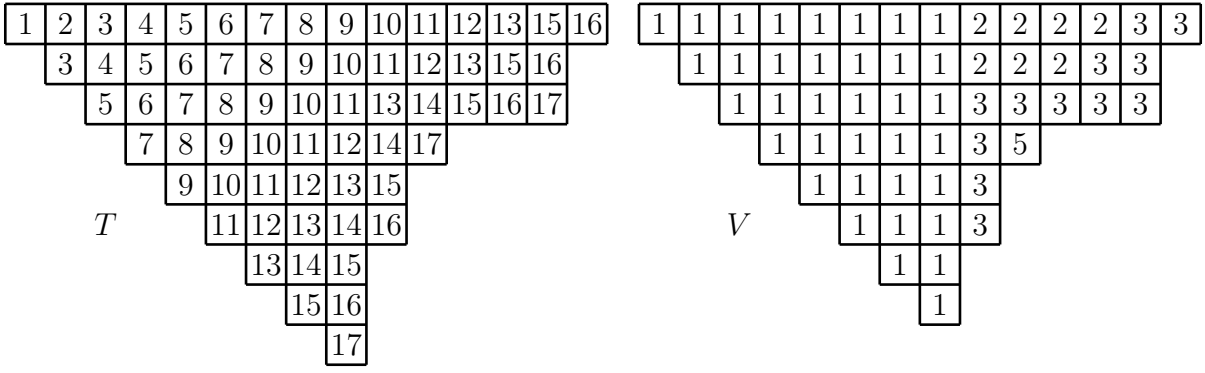
\begin{figure}[ht]
\begin{center}
\begin{tabular}{cc}
\psset{unit=0.21in, linewidth=1pt, radius=2pt, labelsep = 2pt}
\begin{pspicture}(1,0)(16,-9)
\psline[]{-}(1,0)(16,0)
\psline[]{-}(1,-1)(16,-1)
\psline[]{-}(2,-2)(15,-2)
\psline[]{-}(3,-3)(15,-3)
\psline[]{-}(4,-4)(12,-4)
\psline[]{-}(5,-5)(11,-5)
\psline[]{-}(6,-6)(11,-6)
\psline[]{-}(7,-7)(10,-7)
\psline[]{-}(8,-8)(10,-8)
\psline[]{-}(9,-9)(10,-9)
\psline[]{-}(1,0)(1,-1)
\psline[]{-}(2,0)(2,-2)
\psline[]{-}(3,0)(3,-3)
\psline[]{-}(4,0)(4,-4)
\psline[]{-}(5,0)(5,-5)
\psline[]{-}(6,0)(6,-6)
\psline[]{-}(7,0)(7,-7)
\psline[]{-}(8,0)(8,-8)
\psline[]{-}(9,0)(9,-9)
\psline[]{-}(10,0)(10,-9)
\psline[]{-}(11,0)(11,-6)
\psline[]{-}(12,0)(12,-4)
\psline[]{-}(13,0)(13,-3)
\psline[]{-}(14,0)(14,-3)
\psline[]{-}(15,0)(15,-3)
\psline[]{-}(16,0)(16,-1)
\rput(1.5,-0.5){$1$}
\rput(2.5,-0.5){$2$}
\rput(3.5,-0.5){$3$}
\rput(4.5,-0.5){$4$}
\rput(5.5,-0.5){$5$}
\rput(6.5,-0.5){$6$}
\rput(7.5,-0.5){$7$}
\rput(8.5,-0.5){$8$}
\rput(9.5,-0.5){$9$}
\rput(10.5,-0.5){$10$}
\rput(11.5,-0.5){$11$}
\rput(12.5,-0.5){$12$}
\rput(13.5,-0.5){$13$}
\rput(14.5,-0.5){$15$}
\rput(15.5,-0.5){$16$}
\rput(2.5,-1.5){$3$}
\rput(3.5,-1.5){$4$}
\rput(4.5,-1.5){$5$}
\rput(5.5,-1.5){$6$}
\rput(6.5,-1.5){$7$}
\rput(7.5,-1.5){$8$}
\rput(8.5,-1.5){$9$}
\rput(9.5,-1.5){$10$}
\rput(10.5,-1.5){$11$}
\rput(11.5,-1.5){$12$}
\rput(12.5,-1.5){$13$}
\rput(13.5,-1.5){$15$}
\rput(14.5,-1.5){$16$}
\rput(3.5,-2.5){$5$}
\rput(4.5,-2.5){$6$}
\rput(5.5,-2.5){$7$}
\rput(6.5,-2.5){$8$}
\rput(7.5,-2.5){$9$}
\rput(8.5,-2.5){$10$}
\rput(9.5,-2.5){$11$}
\rput(10.5,-2.5){$13$}
\rput(11.5,-2.5){$14$}
\rput(12.5,-2.5){$15$}
\rput(13.5,-2.5){$16$}
\rput(14.5,-2.5){$17$}
\rput(4.5,-3.5){$7$}
\rput(5.5,-3.5){$8$}
\rput(6.5,-3.5){$9$}
\rput(7.5,-3.5){$10$}
\rput(8.5,-3.5){$11$}
\rput(9.5,-3.5){$12$}
\rput(10.5,-3.5){$14$}
\rput(11.5,-3.5){$17$}
\rput(5.5,-4.5){$9$}
\rput(6.5,-4.5){$10$}
\rput(7.5,-4.5){$11$}
\rput(8.5,-4.5){$12$}
\rput(9.5,-4.5){$13$}
\rput(10.5,-4.5){$15$}
\rput(6.5,-5.5){$11$}
\rput(7.5,-5.5){$12$}
\rput(8.5,-5.5){$13$}
\rput(9.5,-5.5){$14$}
\rput(10.5,-5.5){$16$}
\rput(7.5,-6.5){$13$}
\rput(8.5,-6.5){$14$}
\rput(9.5,-6.5){$15$}
\rput(8.5,-7.5){$15$}
\rput(9.5,-7.5){$16$}
\rput(9.5,-8.5){$17$}
\rput(3.5,-5.5){$T$}
\end{pspicture} &
\psset{unit=0.21in, linewidth=1pt, radius=2pt, labelsep = 2pt}
\begin{pspicture}(2,0)(16,-9)
\psline[]{-}(2,0)(16,0)
\psline[]{-}(2,-1)(16,-1)
\psline[]{-}(3,-2)(15,-2)
\psline[]{-}(4,-3)(15,-3)
\psline[]{-}(5,-4)(12,-4)
\psline[]{-}(6,-5)(11,-5)
\psline[]{-}(7,-6)(11,-6)
\psline[]{-}(8,-7)(10,-7)
\psline[]{-}(9,-8)(10,-8)
\psline[]{-}(2,0)(2,-1)
\psline[]{-}(3,0)(3,-2)
\psline[]{-}(4,0)(4,-3)
\psline[]{-}(5,0)(5,-4)
\psline[]{-}(6,0)(6,-5)
\psline[]{-}(7,0)(7,-6)
\psline[]{-}(8,0)(8,-7)
\psline[]{-}(9,0)(9,-8)
\psline[]{-}(10,0)(10,-8)
\psline[]{-}(11,0)(11,-6)
\psline[]{-}(12,0)(12,-4)
\psline[]{-}(13,0)(13,-3)
\psline[]{-}(14,0)(14,-3)
\psline[]{-}(15,0)(15,-3)
\psline[]{-}(16,0)(16,-1)
\rput(2.5,-0.5){$1$}
\rput(3.5,-0.5){$1$}
\rput(4.5,-0.5){$1$}
\rput(5.5,-0.5){$1$}
\rput(6.5,-0.5){$1$}
\rput(7.5,-0.5){$1$}
\rput(8.5,-0.5){$1$}
\rput(9.5,-0.5){$1$}
\rput(10.5,-0.5){$2$}
\rput(11.5,-0.5){$2$}
\rput(12.5,-0.5){$2$}
\rput(13.5,-0.5){$2$}
\rput(14.5,-0.5){$3$}
\rput(15.5,-0.5){$3$}
\rput(3.5,-1.5){$1$}
\rput(4.5,-1.5){$1$}
\rput(5.5,-1.5){$1$}
\rput(6.5,-1.5){$1$}
\rput(7.5,-1.5){$1$}
\rput(8.5,-1.5){$1$}
\rput(9.5,-1.5){$1$}
\rput(10.5,-1.5){$2$}
\rput(11.5,-1.5){$2$}
\rput(12.5,-1.5){$2$}
\rput(13.5,-1.5){$3$}
\rput(14.5,-1.5){$3$}
\rput(4.5,-2.5){$1$}
\rput(5.5,-2.5){$1$}
\rput(6.5,-2.5){$1$}
\rput(7.5,-2.5){$1$}
\rput(8.5,-2.5){$1$}
\rput(9.5,-2.5){$1$}
\rput(10.5,-2.5){$3$}
\rput(11.5,-2.5){$3$}
\rput(12.5,-2.5){$3$}
\rput(13.5,-2.5){$3$}
\rput(14.5,-2.5){$3$}
\rput(5.5,-3.5){$1$}
\rput(6.5,-3.5){$1$}
\rput(7.5,-3.5){$1$}
\rput(8.5,-3.5){$1$}
\rput(9.5,-3.5){$1$}
\rput(10.5,-3.5){$3$}
\rput(11.5,-3.5){$5$}
\rput(6.5,-4.5){$1$}
\rput(7.5,-4.5){$1$}
\rput(8.5,-4.5){$1$}
\rput(9.5,-4.5){$1$}
\rput(10.5,-4.5){$3$}
\rput(7.5,-5.5){$1$}
\rput(8.5,-5.5){$1$}
\rput(9.5,-5.5){$1$}
\rput(10.5,-5.5){$3$}
\rput(8.5,-6.5){$1$}
\rput(9.5,-6.5){$1$}
\rput(9.5,-7.5){$1$}
\rput(4.5,-5.5){$V$}
\end{pspicture}
\end{tabular}
\end{center}
\caption{The bijection between the shifted tableaux $T$ and $V$}\label{fig:thetafig}
\end{figure}

An example of the bijection presented in the proof of the previous proposition is illustrated in Figure~\ref{fig:thetafig}.

\begin{proof}[\bf Proof of Proposition \ref{prop:fdyckprefix}] Set $du P = (k_i)_{i \in [m+2]}$, where $m = |P|_d$, $n = |d u P| = |P| + 2$,
$\lambda(duP) = (\lambda_1, \lambda_2, \ldots, \lambda_{m+1})$, where $\lambda_i = n - i - k_i + 1$, $i \in [m+1]$, and $h = \max \{ h_l : l \in [0,k] \}$, where $h_l = \hv(a_l)$, $l \in [0,k]$.

We first note that formula \eqref{eq:fdyckprefix1} is equivalent to the formula
\begin{equation}\label{eq:fdyckprefixequiv} 
\begin{aligned}
& \textrm{ \# increasing shifted tableaux $T = (t_{ij})$ of shape $\lambda(d uP)$ with $\max(T) = n$} \\ 
= & \textrm{ \# multichains $(W_r)$, $r \in [0, h + k + 1]$,} 
\end{aligned}
\end{equation}
such that $W_0 = P$ and
\begin{enumerate}[label=\roman{*})]

\item $W_r = w_{r0} u w_{r1} \cdots u w_{rk}$, for $r \in [0,h]$, where $w_{rl} = a_l$, for $r \in [0, h - h_l]$ and the multichain $(w_{rl})_{r \in [h-h_l,h]}$ is of type $\V$, for every $l \in [0,k]$.

\item The path $w_{h0} u w_{h1} \cdots u w_{h, h+k-r}$ is a prefix of $W_r$, for $r \in [h, h+k]$.

\end{enumerate}

Indeed, since $\delta(d u P) = n$, from Corollary~\ref{cor:basic}, it follows that the left-hand sides of formulas \eqref{eq:fdyckprefix1} and \eqref{eq:fdyckprefixequiv} are equal. On the other hand, given a sequence $(s_l)_{l \in [0,k]}$ such that $a_l \le s_l$, and a multichain $(w_{rl})$ of type $\V$ from $a_l$ to $s_l$, $r \in [h-h_l, h]$, we set $w_{rl} = a_l$ for every $r \in [0, h-h_l]$, and we define $W_r = w_{r0} u w_{r1} \cdots u w_{rk}$. Furthermore, given a sequence $V_i$, $i \in [k+1]$, with the properties stated in formula \eqref{eq:fdyckprefix1} and a path $V \ge s_0 V_1$, we define $W_r$, for $r \in [h, h+k]$, as $W_r = w_{h0} u w_{h1} \cdots u w_{h, h+k-r} V_{h+k+1-r}$ and $W_{h+k+1} = V \ge W_{h+k}$. Clearly, $(W_r)$, $r \in [0,h+k+1]$, is a multichain satisfying the properties of formula \eqref{eq:fdyckprefixequiv}. Since the converse can be easily verified, we deduce that the right-hand sides of formulas \eqref{eq:fdyckprefix1} and \eqref{eq:fdyckprefixequiv} are equal, so that the two formulas are equivalent.

We will prove formula \eqref{eq:fdyckprefixequiv} bijectively. We can easily check that each cell $(i_l, j_l)$ of $F(duP)$, where $i_l = 1 + \frac{1}{2} \sum\limits_{\nu=0}^l |a_\nu|$ and $j_l = n-1-l- \frac{1}{2} \sum\limits_{\nu=0}^l |a_{\nu}|$, $l \in [-1,k]$, corresponds to the last point $A_l$ of the component $a_l$ of $d u P$ for $l \ge 0$, and $(i_{-1}, j_{-1})$ corresponds to the second point $A_{-1}$ of $d u P$. Clearly, the sequence $(j_l)$, $l \in [-1,k+1]$, with $j_{k+1} = 1$, is a decreasing sequence from $n$ to $1$.

Given an increasing shifted tableau $T = (t_{ij})$ of shape $\lambda(duP)$, with $\max(T) = n$, we define $F_r$, $r \in [0, h + k + 1]$, to be the subdiagram of $F(duP)$ consisting of all cells $(i,j)$ with 
\[ t_{ij} \le i + j + b_{jr}, \]
where
\[ b_{jr} =  \min \{h -r + k, l\} + (h - r) [r \le h],\]
and
\begin{equation}\label{eq:lcases}
l = \begin{cases} 
k, & \textrm{ if $j \in [1,j_k)$}; \\ 
\nu, & \textrm{ if $j \in [j_\nu, j_{\nu-1})$, $\nu \in [0,k]$}; \\
-1, & \textrm{ if $j = n$}.
\end{cases}
\end{equation}

Obviously, since $\lambda_1 = n = \max(T)$, by the increasing property of $T$ it follows that $t_{1j} = j$, for $j \in [n]$, so that each $F_r$ contains the first row of $F(duP)$. Furthermore, we can easily check that $F_r$ is a shifted diagram, and $F_{r+1}$ is a subdiagram of $F_r$, for each $r \in [0,h+k]$. It clearly follows that the path corresponding to $F_r$, $r \in [0,h+k+1]$, is of the form $d u W_r$, where $P \le W_r$, and that $(W_r)$, $r \in [0,h+k+1]$, is a multichain (see Remark \ref{rem:mainremark}\ref{rem:r1}). We will show that this multichain satisfies the properties of formula \eqref{eq:fdyckprefixequiv}.

Firstly, we note that for every $r \in [0,h]$ and $l \in [0,k]$, the path $d u W_r$ passes from all points of $duP$ at height $\le h - r + l$, lying on $a_l$. Indeed, since any such point $A(n + i -j, n - i - j)$ is path-connected with the points $A_l$ and $A_{l-1}$, we have that $j \in [j_l, j_{l-1})$ (see Remark \ref{rem:mainremark}\ref{rem:r4}) and, since $n - i - j \le h - r + l$, we obtain that 
\[ t_{ij} \le n \le i + j + l + h - r = i + j + b_{jr},\] which gives that $(i, j) \in F_r$; so, by Remark \ref{rem:mainremark}\ref{rem:r2}, $d u W_r$ passes from $A$.
In particular, $W_r$, $r \in [0, h]$, passes from all points of $P$ at height $l$ lying on $a_l$  (including the endpoints of $a_l$), so that $W_r$ is uniquely decomposed in the form 
\[ W_r = w_{r0} u w_{r1} \cdots u w_{rk}, \] 
where $w_{rl} \in \mathcal{D}$, $a_l \le w_{rl}$ and $w_{rl}, a_l$ have the same low valleys for every $l \in [0,k]$. Furthermore, for each $l \in [0,k]$ and $r \in [0, h - h_l]$, the component $w_{rl}$ of $W_r$ passes from all points of $P$ at height $\le h_l + l$ lying on $a_l$, so that the paths $w_{rl}$ and $a_l$ have the same valleys; hence, $w_{rl} = a_l$. In particular, $w_{0l} = a_l$ for every $l \in [0,k]$, i.e., $W_0 = P$.
Moreover, since for each $l \in [0,k]$ and for every $r \in [h - h_{l}, h]$ the component $w_{rl}$ passes from all points of $duP$ at height $\le h - r + l$ lying on  $a_l$, the paths $w_{rl}$ and $a_l$ have the same valleys for every height $\le h - r$ (considering that they both start from the origin), so that the multichain $(w_{rl})_{r \in [h-h_l, h]}$ is of type $\V$.

In addition, we will show that for every $r \in [h, h+k]$ the path $w_{h0} u w_{h1} \cdots u w_{h,h+k-r}$ is a prefix of $W_r$.  Indeed, let $A(n+i-j, n-i-j)$ be a point in the component $w_{hl}$, for some $l \in [0,h+k-r]$; then, since $j \in [j_l, j_{l-1})$ and $(i,j) \in F_h$, we have that
\[ t_{ij} \le i + j + l 
= i + j + b_{jr}, \]
so that $(i,j) \in F_r$, and consequently, according to Remark \ref{rem:mainremark}\ref{rem:r2}, $A$ is a point of $W_r$.

Thus, the multichain $(W_r)$, satisfies the properties of formula \eqref{eq:fdyckprefixequiv}.

We finally note that since 
\[ 0 \le b_{jr} - b_{j,r+1} \le 1, \]
we can easily check that 
\begin{equation}\label{eq:tijcond} t_{ij} = i + j + b_{jr}, 
\end{equation}
for every cell $(i, j) \in F_r \setminus F_{r+1}$, $r \in [0, h+k+1]$, where $l$ is given by equality \eqref{eq:lcases} and $F_{h+k+2}$ is the empty diagram. 

For the converse, given a multichain $(W_r)$, $r \in [0, h+k+1]$, satisfying the properties of formula \eqref{eq:fdyckprefixequiv}, in view of the above observation, we define $T = (t_{ij})$ by formula \eqref{eq:tijcond}, where $F_r = F(du W_r)$. It is enough to show that $T$ is an increasing shifted tableau of shape $\lambda (d u P)$ with $\max(T) = n$. We consider two cases:

Firstly, assume that $(i,j) \in F_{h+k+1}$; then, it follows immediately from \eqref{eq:tijcond} that  
$t_{ij} = i + j - 1$ and $(i-1,j)$, $(i,j-1) \in F_{h+k+1}$, so that $t_{i-1, j} = t_{i, j-1} = i + j - 2 < t_{ij}$. Furthermore, in this case we have that $t_{ij} = i + j - 1 \le \lambda_i + 2i - 2 = n + i - k_i -1 \le n$, since $P$ is a Dyck prefix.

Secondly, assume that $(i,j) \in F_r \setminus F_{r+1}$ for some $r \in [0,h+k]$, so that the point $A(n + i - j, n-i-j)$ lies below $duW_{r+1}$ and weakly above $duW_r$. Let $l \in [0,k]$ such that $n + i_{l-1} - j_{l-1} < n + i - j \le n + i_l - j_l$. We consider the following two subcases:

1. Assume that $r < h$ and let $(i^\prime, j^\prime)$, $(i^{\prime\prime}, j^{\prime\prime})$ be the cells of $F(duP)$ corresponding to the points $A^\prime$, $A^{\prime\prime}$ of $duW_r$, $duW_{r+1}$ respectively, which lie on the line $x = n + i - j$ (see Figure \ref{fig:111}). 

\begin{figure}[ht]
\begin{center}
\psset{unit=1em}
\begin{pspicture}(-2,-1)(25,9)

\psline[](-2,0)(-1,-1)(0,0)(4,0)  
\psline[linecolor=red](-2,0.1)(-1,-0.9)(0,0.1)
\psline[linecolor=blue](-2,0.2)(-1,-0.8)(0,0.2)

\psdots[](0,0)(4,0) \psellipticarc[](2,0)(2,1.5){0}{180} \rput[b](2,0.2){$a_0$}
\psline[linestyle=dashed](4,0)(6,2)
\psline[linestyle=dashed, linecolor=red](4,0.1)(6,2.1)
\psline[linestyle=dashed, linecolor=blue](4,0.2)(6,2.2)
 
\psline[](6,2)(10,2)(11,3) \psline[linecolor=red](10,2.1)(11,3.1) \psline[linecolor=blue](10,2.2)(11,3.2)
\psdots[](6,2)(10,2) \rput[lt](10,2){$A_{\ell-1}$} \psellipticarc[](8,2)(2,1.5){0}{180} \rput[b](8,2.2){$a_{\ell-1}$}
\psline[](11,3)(15,3) \psdots[](11,3)(15,3) \rput[lt](15,3){$A_{\ell}$} 
\psellipticarc[](13,3)(2,1.5){0}{180}  \rput[b](13,3.2){$a_{\ell}$}
\psline[linestyle=dashed](15,3)(17,5) 
\psline[linestyle=dashed, linecolor=red](15,3.1)(17,5.1) 
\psline[linestyle=dashed, linecolor=blue](15,3.2)(17,5.2) 

\psline[](17,5)(21,5) \psdots[](17,5)(21,5) \psellipticarc[](19,5)(2,1.5){0}{180} \rput[b](19,5.2){$a_{k}$}

\psellipticarc[linecolor=red](2,0)(2,2){0}{180}  \psellipticarc[linecolor=red](8,2)(2,2.5){0}{180} \psellipticarc[linecolor=red](13,3)(2,2){0}{180} 
\psellipticarc[linecolor=red](19,5)(2,2.5){0}{180} 

\psellipticarc[linecolor=blue](2,0)(2,4){0}{180}  \psellipticarc[linecolor=blue](8,2)(2,4){0}{180} 
\psellipticarc[linecolor=blue](13,3)(2,5){0}{180} \psellipticarc[linecolor=blue](19,5)(2,4){0}{180} 

\psdots[](12.5,4.9)(12.5,5.9)(12.5,7.8)(13.5,6.8) 
\psline[linestyle=dotted](12.5,4.9)(12.5,5.9)(12.5,7.8)
\psline[linestyle=dotted](12.5,5.9)(13.5,6.8) 
\rput[r](12.3,5.2){$A'$} \rput[r](12.2,6.2){$A$}  \rput[rb](12.7,7.8){$A''$} \rput[lb](13.6,5.9){$B$}

\rput[lt](21.2,5){$W_{0}=P$}
\rput[l](21.2,5.7){\red{$W_{r}$}}
\rput[l](21.2,6.8){\blue{$W_{r+1}$}}

\psdots[](0,0)(4,0)(6,2)(10,2)(11,3)(15,3)(17,5)(21,5)
\end{pspicture}
\end{center}
\caption{The decomposition of the paths $duW_r$, $duW_{r+1}$ for $r < h$}\label{fig:111} 
\end{figure}

Then, since the points $A^\prime$, $A^{\prime\prime}$ are path-connected with the points $A_l$, $A_{l-1}$, we deduce that $j^\prime, j^{\prime\prime} \in [j_l, j_{l-1})$ (see Remark \ref{rem:mainremark}\ref{rem:r4}). Furthermore, since $j^{\prime\prime} < j \le j^\prime$ we deduce that $j \in (j_l, j_{l-1})$, so that, by formula \eqref{eq:tijcond}, we have that 
\[ t_{ij} = i + j + l + h - r.\] 
Clearly, since the point $B$ corresponding to the cell $(i,j-1)$ is weakly below $duW_{r+1}$ and above $duW_r$, we deduce that $(i,j-1) \in F_s \setminus F_{s+1}$ for some $s \ge r$ and since $j-1 \in [j_l, j_{l-1})$, by formula \eqref{eq:tijcond}, we have that 
\[ t_{i,j-1} = i + j - 1 + \min\{h - s + k, l\} + (h-s)[s \le h] \le i + j - 1 + l + h-r < t_{ij}.\] 
Similarly, we can show that $t_{i-1,j} < t_{ij}$. 

Now, since the point $A$ lies below the component $w_{r+1,l}$ and weakly above the component $w_{rl}$, we have that $w_{r+1,l} \ne a_l$, so that $r \ge h - h_l$. Furthermore, since the sequence $(w_{rl})$, $r \in [h-h_l, h]$, is of type $\V$, we have that the paths $w_{rl}$, $w_{r+1,l}$ have the same valleys at every height $\le h_l - (r - h + h_l + 1) = h - r - 1$ (considering that they both start from the origin), i.e., the components $w_{rl}$, $w_{r+1,l}$ have the same valleys at every height $\le h - r - 1 + l$. Thus, since $A$ does not belong to the component $w_{r+1,l}$, its height is greater than $h - r - 1 + l$, i.e., $n - i - j  > h - r - 1 + l$, so that $t_{ij} \le n$.

2. Assume that $r \ge h$. Clearly, by property ii) of \eqref{eq:fdyckprefixequiv}, the paths $duW_r$ and $duW_{r+1}$ coincide for every point up to $A_{h+k-r-1}$, and since $A$ lies below $duW_{r+1}$ and weakly above $duW_r$, the point $A_l$ must lie on the right of $A_{h+k-r-1}$, i.e., $l \ge h + k - r$ (see Figure \ref{fig:222}).

\begin{figure}[ht]
\psset{unit=1em}
\begin{pspicture}(-2,-1)(36,14)

\psline[](-2,0)(-1,-1)(0,0)(4,0)  
\psline[linecolor=green](-2,0.1)(-1,-0.9)(0,0.1)
\psline[linecolor=red](-2,0.2)(-1,-0.8)(0,0.2)
\psline[linecolor=blue](-2,0.3)(-1,-0.7)(0,0.3)

\psdots[](0,0)(4,0) \psellipticarc[](2,0)(2,1.5){0}{180} \rput[b](2,0.2){$a_0$}
\psline[linestyle=dashed](4,0)(6,2)
\psline[linestyle=dashed, linecolor=green](4,0.1)(6,2.1)
\psline[linestyle=dashed, linecolor=red](4,0.2)(6,2.2)
\psline[linestyle=dashed, linecolor=blue](4,0.3)(6,2.3)
\psline[](6,2)(10,2)(11,3) 
\psline[linecolor=green](10,2.1)(11,3.1) 
\psline[linecolor=red](10,2.2)(11,3.2) 
\psline[linecolor=blue](10,2.3)(11,3.3) 

\psdots[](6,2)(10,2) \rput[lt](10,2){$A_{h+k-r-1}$} 
\psellipticarc[](8,2)(2,1.5){0}{180} \rput[b](8,2.2){$a_{h+k\!-\!r\!-\! 1}$}

\psline[](11,3)(15,3)
 \psellipticarc[](13,3)(2,1.5){0}{180} \rput[b](13,3.2){$a_{h+k-r}$}
 \psdots[](15,3) \rput[lt](15,3){$A_{h+k-r}$}

\psline[linestyle=dashed](15,3)(17,5)
 
\psline[](17,5)(21,5)(22,6) 
\psline[linecolor=green](21,5.1)(22,6.1) 
\psdots[](17,5)(21,5) \rput[lt](21,5){$A_{\ell-1}$} 
\psellipticarc[](19,5)(2,1.5){0}{180} \rput[b](19,5.2){$a_{\ell-1}$}
\psline[](22,6)(26,6) \psdots[](22,6)(26,6) \rput[lt](26,6){$A_{\ell}$} 
\psellipticarc[](24,6)(2,1.5){0}{180}  \rput[b](24,6.2){$a_{\ell}$}
\psline[linestyle=dashed](26,6)(28,8) 
\psline[linestyle=dashed, linecolor=green](26,6.1)(28,8.1) 
\psline[](28,8)(32,8) \psdots[](28,8)(32,8) \psellipticarc[](30,8)(2,1.5){0}{180} \rput[b](30,8.2){$a_{k}$}

\psellipticarc[linecolor=green](2,0)(2,2){0}{180}  \psellipticarc[linecolor=green](8,2)(2,2){0}{180} 
\psellipticarc[linecolor=green](13,3)(2,2){0}{180}
\psline[linestyle=dashed, linecolor = green](15,3.1)(17,5.1)
\psellipticarc[linecolor=green](19,5)(2,2){0}{180} 
\psellipticarc[linecolor=green](24,6)(2,2){0}{180} 
\psellipticarc[linecolor=green](30,8)(2,2){0}{180} 

\psellipticarc[linecolor=red](2,0)(2,3){0}{180}  \psellipticarc[linecolor=red](8,2)(2,3){0}{180} 
\psellipticarc[linecolor=red](13,3)(2,3){0}{180}  
\psline[linestyle=dashed, linecolor = red](15,3.2)(17,5.2)
\psellipticarc[linecolor=red](28,5)(11,6){54}{180} 

\psellipticarc[linecolor=blue](2,0)(2,4){0}{180}  \psellipticarc[linecolor=blue](8,2)(2,4){0}{180} 
\psellipticarc[linecolor=blue](25,3)(14,11){53}{180}  

\psdots[](24,8)(24,12)(25,13) \psline[linestyle=dotted](24,8)(24,12)(25,13)
\rput[r](23.8,8.3){$A'$} \rput[r](23.5,12.2){$A$}   \rput[lb](25.2,12.5){$B$}

\rput[l](32.2,7){$W_0=P$}
\rput[l](32.2,9){\green{$W_h$}}
\rput[l](32.2,10.5){\red{$W_r$}}
\rput[l](32.2,12.5){\blue{$W_{r+1}$}}

\psdots[](0,0)(4,0)(6,2)(10,2)(11,3)(15,3)(17,5)(21,5)(22,6)(26,6)(28,8)(32,8)
\end{pspicture}
\caption{The decomposition of the paths $duW_r$, $duW_{r+1}$ for $r \ge h$ }\label{fig:222} 
\end{figure}

Furthermore, if $(i^\prime,j^\prime)$ is the cell corresponding to the point $A^{\prime}$ which lies on the path $duW_h$ and on the line $x = n + i - j$, then $j^\prime > j$ and $j^\prime \in [j_l, j_{l-1})$, so that $j \in [j_{\nu}, j_{\nu-1})$, for some $\nu \in [l,k+1]$. Then, by formula \eqref{eq:tijcond} we have that 
\[ 
t_{ij} = i + j + \min \{h -r + k, \nu\} + (h-r)[r \le h] = i + j + h - r +k. 
\]
Clearly, as in the previous subcase, we have $(i,j-1) \in F_s \setminus F_{s+1}$ for some $s \ge r$ and so, by the above formula, we deduce that
\[ t_{i,j-1} = i + j - 1 + h - s + k < t_{ij}.\] 
Similarly, we obtain that $t_{i-1,j} < t_{ij}$.  

Finally, since the height of $A$ is greater than or equal to $l$, we have that $n - i - j \ge l \ge h - r + k$, so that $t_{ij} \le n$.
\end{proof}

An example of the bijection used for the proof of formula \eqref{eq:fdyckprefixequiv} is illustrated in Figure~\ref{fig:wrmultichain}, where the cells $(i, j) \in F(duW_r) \setminus F(duW_{r+1})$ of the tableau are colored with the same color as $duW_r$.

\newpage

\begin{figure}[ht]
\begin{center}
\psset{unit=.65em, linewidth=1pt, radius=2pt, labelsep = 2pt}
\begin{tabular}{c}
\begin{tabular}{c}
\psset{unit=0.18in, linewidth=1pt, radius=2pt, labelsep = 2pt}
{\scriptsize
\begin{pspicture}(2.75,-20)(39,0)
\psline[linestyle=dotted]{-}(38,-1)(34,-5)
\psline[linestyle=dotted]{-}(33,-5)(29,-9)
\psline[linestyle=dotted]{-}(28,-9)(21,-16)
\psline[linestyle=dotted]{-}(20,-16)(18,-18)
\psline[]{-}(0,0)(39,-0)
\psline[]{-}(0,-1)(39,-1)
\psline[]{-}(1,-2)(36,-2)
\psline[]{-}(2,-3)(34,-3)
\psline[]{-}(3,-4)(34,-4)
\psline[]{-}(4,-5)(34,-5)
\psline[]{-}(5,-6)(32,-6)
\psline[]{-}(6,-7)(30,-7)
\psline[]{-}(7,-8)(29,-8)
\psline[]{-}(8,-9)(29,-9)
\psline[]{-}(9,-10)(27,-10)
\psline[]{-}(10,-11)(26,-11)
\psline[]{-}(11,-12)(24,-12)
\psline[]{-}(12,-13)(22,-13)
\psline[]{-}(13,-14)(21,-14)
\psline[]{-}(14,-15)(21,-15)
\psline[]{-}(15,-16)(21,-16)
\psline[]{-}(16,-17)(19,-17)
\psline[]{-}(17,-18)(18,-18)
\psline[]{-}(0,0)(0,-1)
\psline[]{-}(1,0)(1,-2)
\psline[]{-}(2,0)(2,-3)
\psline[]{-}(3,0)(3,-4)
\psline[]{-}(4,0)(4,-5)
\psline[]{-}(5,0)(5,-6)
\psline[]{-}(6,0)(6,-7)
\psline[]{-}(7,0)(7,-8)
\psline[]{-}(8,0)(8,-9)
\psline[]{-}(9,0)(9,-10)
\psline[]{-}(10,0)(10,-11)
\psline[]{-}(11,0)(11,-12)
\psline[]{-}(12,0)(12,-13)
\psline[]{-}(13,0)(13,-14)
\psline[]{-}(14,0)(14,-15)
\psline[]{-}(15,0)(15,-16)
\psline[]{-}(16,0)(16,-17)
\psline[]{-}(17,0)(17,-18)
\psline[]{-}(18,0)(18,-18)
\psline[]{-}(19,0)(19,-17)
\psline[]{-}(20,0)(20,-16)
\psline[]{-}(21,0)(21,-16)
\psline[]{-}(22,0)(22,-13)
\psline[]{-}(23,0)(23,-12)
\psline[]{-}(24,0)(24,-12)
\psline[]{-}(25,0)(25,-11)
\psline[]{-}(26,0)(26,-11)
\psline[]{-}(27,0)(27,-10)
\psline[]{-}(28,0)(28,-9)
\psline[]{-}(29,0)(29,-9)
\psline[]{-}(30,0)(30,-7)
\psline[]{-}(31,0)(31,-6)
\psline[]{-}(32,0)(32,-6)
\psline[]{-}(33,0)(33,-5)
\psline[]{-}(34,0)(34,-5)
\psline[]{-}(35,0)(35,-2)
\psline[]{-}(36,0)(36,-2)
\psline[]{-}(37,0)(37,-1)
\psline[]{-}(38,0)(38,-1)
\psline[]{-}(39,0)(39,-1)
\rput(0.5,-0.5){$\blue 1$}
\rput(1.5,-0.5){$\blue 2$}
\rput(2.5,-0.5){$\blue 3$}
\rput(3.5,-0.5){$\blue 4$}
\rput(4.5,-0.5){$\blue 5$}
\rput(5.5,-0.5){$\blue 6$}
\rput(6.5,-0.5){$\blue 7$}
\rput(7.5,-0.5){$\blue 8$}
\rput(8.5,-0.5){$\blue 9$}
\rput(9.5,-0.5){$\blue 10$}
\rput(10.5,-0.5){$\blue 11$}
\rput(11.5,-0.5){$\blue 12$}
\rput(12.5,-0.5){$\blue 13$}
\rput(13.5,-0.5){$\blue 14$}
\rput(14.5,-0.5){$\blue 15$}
\rput(15.5,-0.5){$\blue 16$}
\rput(16.5,-0.5){$\blue 17$}
\rput(17.5,-0.5){$\blue 18$}
\rput(18.5,-0.5){$\blue 19$}
\rput(19.5,-0.5){$\blue 20$}
\rput(20.5,-0.5){$\blue 21$}
\rput(21.5,-0.5){$\blue 22$}
\rput(22.5,-0.5){$\blue 23$}
\rput(23.5,-0.5){$\blue 24$}
\rput(24.5,-0.5){$\blue 25$}
\rput(25.5,-0.5){$\blue 26$}
\rput(26.5,-0.5){$\blue 27$}
\rput(27.5,-0.5){$\blue 28$}
\rput(28.5,-0.5){$\blue 29$}
\rput(29.5,-0.5){$\blue 30$}
\rput(30.5,-0.5){$\blue 31$}
\rput(31.5,-0.5){$\blue 32$}
\rput(32.5,-0.5){$\blue 33$}
\rput(33.5,-0.5){$\blue 34$}
\rput(34.5,-0.5){$\blue 35$}
\rput(35.5,-0.5){$\blue 36$}
\rput(36.5,-0.5){$\blue 37$}
\rput(37.5,-0.5){$\blue 38$}
\rput(38.5,-0.5){$\blue 39$}
\rput(1.5,-1.5){$\blue 3$}
\rput(2.5,-1.5){$\blue 4$}
\rput(3.5,-1.5){$\blue 5$}
\rput(4.5,-1.5){$\blue 6$}
\rput(5.5,-1.5){$\blue 7$}
\rput(6.5,-1.5){$\blue 8$}
\rput(7.5,-1.5){$\blue 9$}
\rput(8.5,-1.5){$\blue 10$}
\rput(9.5,-1.5){$\blue 11$}
\rput(10.5,-1.5){$\blue 12$}
\rput(11.5,-1.5){$\blue 13$}
\rput(12.5,-1.5){$\blue 14$}
\rput(13.5,-1.5){$\blue 15$}
\rput(14.5,-1.5){$\blue 16$}
\rput(15.5,-1.5){$\blue 17$}
\rput(16.5,-1.5){$\blue 18$}
\rput(17.5,-1.5){$\blue 19$}
\rput(18.5,-1.5){$\blue 20$}
\rput(19.5,-1.5){$\blue 21$}
\rput(20.5,-1.5){$\blue 22$}
\rput(21.5,-1.5){$\blue 23$}
\rput(22.5,-1.5){$\blue 24$}
\rput(23.5,-1.5){$\blue 25$}
\rput(24.5,-1.5){$\blue 26$}
\rput(25.5,-1.5){$\blue 27$}
\rput(26.5,-1.5){$\blue 28$}
\rput(27.5,-1.5){$\blue 29$}
\rput(28.5,-1.5){$\blue 30$}
\rput(29.5,-1.5){$\blue 31$}
\rput(30.5,-1.5){$\blue 32$}
\rput(31.5,-1.5){$\blue 33$}
\rput(32.5,-1.5){${\red 35}$}
\rput(33.5,-1.5){${\red 36}$}
\rput(34.5,-1.5){${\yellow 38}$}
\rput(35.5,-1.5){${\yellow 39}$}
\rput(2.5,-2.5){$\blue 5$}
\rput(3.5,-2.5){$\blue 6$}
\rput(4.5,-2.5){$\blue 7$}
\rput(5.5,-2.5){$\blue 8$}
\rput(6.5,-2.5){$\blue 9$}
\rput(7.5,-2.5){$\blue 10$}
\rput(8.5,-2.5){$\blue 11$}
\rput(9.5,-2.5){$\blue 12$}
\rput(10.5,-2.5){$\blue 13$}
\rput(11.5,-2.5){$\blue 14$}
\rput(12.5,-2.5){$\blue 15$}
\rput(13.5,-2.5){$\blue 16$}
\rput(14.5,-2.5){$\blue 17$}
\rput(15.5,-2.5){$\blue 18$}
\rput(16.5,-2.5){$\blue 19$}
\rput(17.5,-2.5){$\blue 20$}
\rput(18.5,-2.5){$\blue 21$}
\rput(19.5,-2.5){$\blue 22$}
\rput(20.5,-2.5){$\blue 23$}
\rput(21.5,-2.5){$\blue 24$}
\rput(22.5,-2.5){$\blue 25$}
\rput(23.5,-2.5){$\blue 26$}
\rput(24.5,-2.5){$\blue 27$}
\rput(25.5,-2.5){$\blue 28$}
\rput(26.5,-2.5){$\blue 29$}
\rput(27.5,-2.5){$\blue 30$}
\rput(28.5,-2.5){$\blue 31$}
\rput(29.5,-2.5){$\blue 32$}
\rput(30.5,-2.5){$\blue 33$}
\rput(31.5,-2.5){${\red 35}$}
\rput(32.5,-2.5){${\red 36} $}
\rput(33.5,-2.5){${\red 37}$}
\rput(3.5,-3.5){$\blue 7$}
\rput(4.5,-3.5){$\blue 8$}
\rput(5.5,-3.5){$\blue 9$}
\rput(6.5,-3.5){$\blue 10$}
\rput(7.5,-3.5){$\blue 11$}
\rput(8.5,-3.5){$\blue 12$}
\rput(9.5,-3.5){$\blue 13$}
\rput(10.5,-3.5){$\blue 14$}
\rput(11.5,-3.5){$\blue 15$}
\rput(12.5,-3.5){$\blue 16$}
\rput(13.5,-3.5){$\blue 17$}
\rput(14.5,-3.5){$\blue 18$}
\rput(15.5,-3.5){$\blue 19$}
\rput(16.5,-3.5){$\blue 20$}
\rput(17.5,-3.5){$\blue 21$}
\rput(18.5,-3.5){$\blue 22$}
\rput(19.5,-3.5){$\blue 23$}
\rput(20.5,-3.5){$\blue 24$}
\rput(21.5,-3.5){$\blue 25$}
\rput(22.5,-3.5){$\blue 26$}
\rput(23.5,-3.5){$\blue 27$}
\rput(24.5,-3.5){$\blue 28$}
\rput(25.5,-3.5){${\red 30}$}
\rput(26.5,-3.5){${\red 31}$}
\rput(27.5,-3.5){${\red 32}$}
\rput(28.5,-3.5){${\red 33}$}
\rput(29.5,-3.5){${\red 34}$}
\rput(30.5,-3.5){${\red 35}$}
\rput(31.5,-3.5){${\red 36}$}
\rput(32.5,-3.5){${\red 37}$}
\rput(33.5,-3.5){${\red 38}$}
\rput(4.5,-4.5){$\blue 9$}
\rput(5.5,-4.5){$\blue 10$}
\rput(6.5,-4.5){$\blue 11$}
\rput(7.5,-4.5){$\blue 12$}
\rput(8.5,-4.5){$\blue 13$}
\rput(9.5,-4.5){$\blue 14$}
\rput(10.5,-4.5){$\blue 15$}
\rput(11.5,-4.5){$\blue 16$}
\rput(12.5,-4.5){$\blue 17$}
\rput(13.5,-4.5){$\blue 18$}
\rput(14.5,-4.5){$\blue 19$}
\rput(15.5,-4.5){$\blue 20$}
\rput(16.5,-4.5){$\blue 21$}
\rput(17.5,-4.5){$\blue 22$}
\rput(18.5,-4.5){$\blue 23$}
\rput(19.5,-4.5){$\blue 24$}
\rput(20.5,-4.5){$\blue 25$}
\rput(21.5,-4.5){$\blue 26$}
\rput(22.5,-4.5){$\blue 27$}
\rput(23.5,-4.5){$\blue 28$}
\rput(24.5,-4.5){$\blue 29$}
\rput(25.5,-4.5){${\red 31}$}
\rput(26.5,-4.5){${\red 32}$}
\rput(27.5,-4.5){${\red 33}$}
\rput(28.5,-4.5){${\red 34}$}
\rput(29.5,-4.5){${\red 35}$}
\rput(30.5,-4.5){${\red 36}$}
\rput(31.5,-4.5){${\red 37}$}
\rput(32.5,-4.5){${\red 38}$}
\rput(33.5,-4.5){${\red 39}$}
\rput(5.5,-5.5){$\blue 11$}
\rput(6.5,-5.5){$\blue 12$}
\rput(7.5,-5.5){$\blue 13$}
\rput(8.5,-5.5){$\blue 14$}
\rput(9.5,-5.5){$\blue 15$}
\rput(10.5,-5.5){$\blue 16$}
\rput(11.5,-5.5){$\blue 17$}
\rput(12.5,-5.5){$\blue 18$}
\rput(13.5,-5.5){$\blue 19$}
\rput(14.5,-5.5){$\blue 20$}
\rput(15.5,-5.5){$\blue 21$}
\rput(16.5,-5.5){$\blue 22$}
\rput(17.5,-5.5){$\blue 23$}
\rput(18.5,-5.5){$\blue 24$}
\rput(19.5,-5.5){$\blue 25$}
\rput(20.5,-5.5){$\blue 26$}
\rput(21.5,-5.5){$\blue 27$}
\rput(22.5,-5.5){$\blue 28$}
\rput(23.5,-5.5){$\blue 29$}
\rput(24.5,-5.5){$\blue 30$}
\rput(25.5,-5.5){${\red 32}$}
\rput(26.5,-5.5){${\red 33}$}
\rput(27.5,-5.5){${\red 34}$}
\rput(28.5,-5.5){${\red 35}$}
\rput(29.5,-5.5){${\green 37}$}
\rput(30.5,-5.5){${\green 38}$}
\rput(31.5,-5.5){${\green 39}$}
\rput(6.5,-6.5){$\blue 13$}
\rput(7.5,-6.5){$\blue 14$}
\rput(8.5,-6.5){$\blue 15$}
\rput(9.5,-6.5){$\blue 16$}
\rput(10.5,-6.5){$\blue 17$}
\rput(11.5,-6.5){$\blue 18$}
\rput(12.5,-6.5){$\blue 19$}
\rput(13.5,-6.5){$\blue 20$}
\rput(14.5,-6.5){$\blue 21$}
\rput(15.5,-6.5){$\blue 22$}
\rput(16.5,-6.5){$\blue 23$}
\rput(17.5,-6.5){$\blue 24$}
\rput(18.5,-6.5){$\blue 25$}
\rput(19.5,-6.5){$\blue 26$}
\rput(20.5,-6.5){$\blue 27$}
\rput(21.5,-6.5){$\blue 28$}
\rput(22.5,-6.5){$\blue 29$}
\rput(23.5,-6.5){$\blue 30$}
\rput(24.5,-6.5){$\blue 31$}
\rput(25.5,-6.5){${\red 33}$}
\rput(26.5,-6.5){${\red 34}$}
\rput(27.5,-6.5){${\red 35}$}
\rput(28.5,-6.5){${\green 37}$}
\rput(29.5,-6.5){${\green 38}$}
\rput(7.5,-7.5){$\blue 15$}
\rput(8.5,-7.5){$\blue 16$}
\rput(9.5,-7.5){$\blue 17$}
\rput(10.5,-7.5){$\blue 18$}
\rput(11.5,-7.5){$\blue 19$}
\rput(12.5,-7.5){$\blue 20$}
\rput(13.5,-7.5){$\blue 21$}
\rput(14.5,-7.5){$\blue 22$}
\rput(15.5,-7.5){$\blue 23$}
\rput(16.5,-7.5){$\blue 24$}
\rput(17.5,-7.5){$\blue 25$}
\rput(18.5,-7.5){$\blue 26$}
\rput(19.5,-7.5){$\blue 27$}
\rput(20.5,-7.5){$\blue 28$}
\rput(21.5,-7.5){$\blue 29$}
\rput(22.5,-7.5){$\blue 30$}
\rput(23.5,-7.5){${\red 32}$}
\rput(24.5,-7.5){${\red 33}$}
\rput(25.5,-7.5){${\red 34}$}
\rput(26.5,-7.5){${\red 35}$}
\rput(27.5,-7.5){${\red 36}$}
\rput(28.5,-7.5){${\green 38}$}
\rput(8.5,-8.5){$\blue 17$}
\rput(9.5,-8.5){$\blue 18$}
\rput(10.5,-8.5){$\blue 19$}
\rput(11.5,-8.5){$\blue 20$}
\rput(12.5,-8.5){$\blue 21$}
\rput(13.5,-8.5){$\blue 22$}
\rput(14.5,-8.5){$\blue 23$}
\rput(15.5,-8.5){$\blue 24$}
\rput(16.5,-8.5){$\blue 25$}
\rput(17.5,-8.5){$\blue 26$}
\rput(18.5,-8.5){$\blue 27$}
\rput(19.5,-8.5){$\blue 28$}
\rput(20.5,-8.5){$\blue 29$}
\rput(21.5,-8.5){$\blue 30$}
\rput(22.5,-8.5){${\red 32}$}
\rput(23.5,-8.5){${\red 33}$}
\rput(24.5,-8.5){${\red 34}$}
\rput(25.5,-8.5){${\red 35}$}
\rput(26.5,-8.5){${\red 36}$}
\rput(27.5,-8.5){${\red 37}$}
\rput(28.5,-8.5){${\green 39}$}
\rput(9.5,-9.5){$\blue 19$}
\rput(10.5,-9.5){$\blue 20$}
\rput(11.5,-9.5){$\blue 21$}
\rput(12.5,-9.5){$\blue 22$}
\rput(13.5,-9.5){$\blue 23$}
\rput(14.5,-9.5){$\blue 24$}
\rput(15.5,-9.5){$\blue 25$}
\rput(16.5,-9.5){$\blue 26$}
\rput(17.5,-9.5){$\blue 27$}
\rput(18.5,-9.5){$\blue 28$}
\rput(19.5,-9.5){$\blue 29$}
\rput(20.5,-9.5){$\blue 30$}
\rput(21.5,-9.5){$\blue 31$}
\rput(22.5,-9.5){${\red 33}$}
\rput(23.5,-9.5){${\red 34}$}
\rput(24.5,-9.5){${\red 35}$}
\rput(25.5,-9.5){${\red 36}$}
\rput(26.5,-9.5){${\red 37}$}
\rput(10.5,-10.5){$\blue 21$}
\rput(11.5,-10.5){$\blue 22$}
\rput(12.5,-10.5){$\blue 23$}
\rput(13.5,-10.5){$\blue 24$}
\rput(14.5,-10.5){$\blue 25$}
\rput(15.5,-10.5){$\blue 26$}
\rput(16.5,-10.5){$\blue 27$}
\rput(17.5,-10.5){$\blue 28$}
\rput(18.5,-10.5){$\blue 29$}
\rput(19.5,-10.5){$\blue 30$}
\rput(20.5,-10.5){$\blue 31$}
\rput(21.5,-10.5){$\blue 32$}
\rput(22.5,-10.5){${\red 34}$}
\rput(23.5,-10.5){${\cyan 37}$}
\rput(24.5,-10.5){${\cyan 38}$}
\rput(25.5,-10.5){${\cyan 39}$}
\rput(11.5,-11.5){$\blue 23$}
\rput(12.5,-11.5){$\blue 24$}
\rput(13.5,-11.5){$\blue 25$}
\rput(14.5,-11.5){$\blue 26$}
\rput(15.5,-11.5){$\blue 27$}
\rput(16.5,-11.5){$\blue 28$}
\rput(17.5,-11.5){$\blue 29$}
\rput(18.5,-11.5){${\red 31}$}
\rput(19.5,-11.5){${\red 32}$}
\rput(20.5,-11.5){${\red 33}$}
\rput(21.5,-11.5){${\green 35}$}
\rput(22.5,-11.5){${\yellow 38}$}
\rput(23.5,-11.5){${\yellow 39}$}
\rput(12.5,-12.5){$\blue 25$}
\rput(13.5,-12.5){$\blue 26$}
\rput(14.5,-12.5){${\red 28}$}
\rput(15.5,-12.5){${\green 30}$}
\rput(16.5,-12.5){${\green 31}$}
\rput(17.5,-12.5){${\green 32}$}
\rput(18.5,-12.5){${\green 33}$}
\rput(19.5,-12.5){${\green 34}$}
\rput(20.5,-12.5){${\green 35}$}
\rput(21.5,-12.5){${\black 39}$}
\rput(13.5,-13.5){${\green 29}$}
\rput(14.5,-13.5){${\green 30}$}
\rput(15.5,-13.5){${\green 31}$}
\rput(16.5,-13.5){${\green 32}$}
\rput(17.5,-13.5){${\green 33}$}
\rput(18.5,-13.5){${\green 34}$}
\rput(19.5,-13.5){${\green 35}$}
\rput(20.5,-13.5){${\green 36}$}
\rput(14.5,-14.5){${\green 31}$}
\rput(15.5,-14.5){${\green 32}$}
\rput(16.5,-14.5){${\green 33}$}
\rput(17.5,-14.5){${\green 34}$}
\rput(18.5,-14.5){${\green 35}$}
\rput(19.5,-14.5){${\green 36}$}
\rput(20.5,-14.5){${\cyan 38}$}
\rput(15.5,-15.5){${\green 33}$}
\rput(16.5,-15.5){${\green 34}$}
\rput(17.5,-15.5){${\cyan 36}$}
\rput(18.5,-15.5){${\cyan 37}$}
\rput(19.5,-15.5){${\cyan 38}$}
\rput(20.5,-15.5){${\cyan 39}$}
\rput(16.5,-16.5){${\cyan 36}$}
\rput(17.5,-16.5){${\cyan 37}$}
\rput(18.5,-16.5){${\yellow 39}$}
\rput(17.5,-17.5){${\magenta 39}$}
\end{pspicture}
}
\end{tabular} \\
\begin{tabular}{c}
\begin{pspicture}(8,-1)(40,12)
\rput(-2.5,0){\footnotesize $duW_6$}
\psline[showpoints=true,linecolor=blue]{-}(0,0)(1,-1)(2,0)(3,1)(4,2)(5,3)(6,4)(7,5)(8,6)(9,5)(10,6)(11,5)(12,6)(13,7)(14,8)(15,9)(16,10)(17,11)(18,10)(19,9)(20,8)(21,7)(22,8)(23,9)(24,8)(25,9)(26,8)(27,7)(28,6)(29,7)(30,8)(31,9)(32,10)(33,9)(34,10)(35,11)(36,12)(37,13)(38,12)(39,13)
\end{pspicture} \\[-1em]
\begin{pspicture}(8,-1)(40,1)
\rput(-2.5,0){\footnotesize $duW_5$}
\psline[showpoints=true,linecolor=red]{-}(0,0)(1,-1)(2,0)(3,1)(4,2)(5,3)(6,4)(7,3)(8,2)(9,1)(10,0)(11,1)(12,2)(13,3)(14,4)(15,5)(16,4)(17,5)(18,4)(19,3)(20,2)(21,3)(22,2)(23,3)(24,4)(25,5)(26,6)(27,5)(28,6)(29,7)(30,6)(31,7)(32,8)(33,9)(34,10)(35,11)(36,12)(37,11)(38,12)(39,13)
\end{pspicture} \\[-1em]
\begin{pspicture}(8,-1)(40,1)
\rput(-2.5,0){\footnotesize $duW_4$}
\psline[showpoints=true,linecolor=green]{-}(0,0)(1,-1)(2,0)(3,1)(4,2)(5,3)(6,4)(7,3)(8,2)(9,1)(10,0)(11,1)(12,2)(13,1)(14,2)(15,3)(16,2)(17,3)(18,2)(19,1)(20,2)(21,3)(22,2)(23,3)(24,4)(25,5)(26,6)(27,5)(28,6)(29,5)(30,6)(31,5)(32,4)(33,5)(34,4)(35,5)(36,6)(37,7)(38,6)(39,7)
\end{pspicture} \\[-1em]
\begin{pspicture}(8,-1)(40,1)
\rput(-2.5,0){\footnotesize $duW_3$}
\psline[showpoints=true,linecolor=cyan]{-}(0,0)(1,-1)(2,0)(3,1)(4,2)(5,3)(6,4)(7,3)(8,2)(9,1)(10,0)(11,1)(12,2)(13,1)(14,2)(15,3)(16,2)(17,3)(18,2)(19,1)(20,2)(21,3)(22,2)(23,3)(24,2)(25,3)(26,4)(27,5)(28,6)(29,5)(30,6)(31,5)(32,4)(33,3)(34,2)(35,3)(36,4)(37,5)(38,4)(39,5)
\end{pspicture} \\[-1em]
\begin{pspicture}(8,-1)(40,1)
\rput(-2.5,0){\footnotesize $duW_2$}
\psline[showpoints=true,linecolor=magenta]{-}(0,0)(1,-1)(2,0)(3,1)(4,2)(5,3)(6,4)(7,3)(8,2)(9,1)(10,0)(11,1)(12,2)(13,1)(14,2)(15,3)(16,2)(17,3)(18,2)(19,1)(20,2)(21,3)(22,2)(23,3)(24,2)(25,3)(26,4)(27,5)(28,6)(29,5)(30,6)(31,5)(32,4)(33,3)(34,2)(35,3)(36,4)(37,5)(38,4)(39,3)
\end{pspicture} \\[-1em]
\begin{pspicture}(8,-1)(40,1)
\rput(-2.5,0){\footnotesize $duW_1$}
\psline[showpoints=true,linecolor=yellow]{-}(0,0)(1,-1)(2,0)(3,1)(4,2)(5,1)(6,2)(7,3)(8,2)(9,1)(10,0)(11,1)(12,2)(13,1)(14,2)(15,3)(16,2)(17,3)(18,2)(19,1)(20,2)(21,3)(22,2)(23,3)(24,2)(25,3)(26,4)(27,3)(28,4)(29,5)(30,6)(31,5)(32,4)(33,3)(34,2)(35,3)(36,4)(37,3)(38,4)(39,3)
\end{pspicture} \\[-1em]
\begin{pspicture}(8,-1)(40,1)
\rput(-2.5,0){\footnotesize $duW_0$}
\psline[showpoints=true,linecolor=black]{-}(0,0)(1,-1)(2,0)(3,1)(4,2)(5,1)(6,2)(7,3)(8,2)(9,1)(10,0)(11,1)(12,2)(13,1)(14,2)(15,3)(16,2)(17,3)(18,2)(19,1)(20,2)(21,3)(22,2)(23,3)(24,2)(25,3)(26,4)(27,3)(28,4)(29,5)(30,4)(31,5)(32,4)(33,3)(34,2)(35,3)(36,4)(37,3)(38,4)(39,3)
\psline[showpoints=true,linestyle=dotted]{-}(2,0)(10,0)
\psline[showpoints=true,linestyle=dotted]{-}(11,1)(19,1)
\psline[showpoints=true,linestyle=dotted]{-}(20,2)(34,2)
\psline[showpoints=true,linestyle=dotted]{-}(35,3)(39,3)
\end{pspicture}
\end{tabular} \\
\end{tabular} 
\end{center}
\caption{An increasing shifted tableau and its corresponding multichain $(duW_r)_{r \in [0,6]}$} \label{fig:wrmultichain}
\end{figure}



\begin{thebibliography}{99}

\bibitem{BP1995}
M. Barnabei and E. Pezzoli, \newblock M\"obius functions, \newblock in J. P. S. Kung ed., {\em Gian-Carlo Rota on combinatorics}, Birkhauser, 1995, 83--104.

\bibitem{CTY2014}
E. Clifford, H. Thomas, and A. Yong, \newblock $K$-theoretic Schubert calculus for $OG(n,2n+1)$ and jeu de taquin for shifted increasing tableaux, \newblock {\em J. Reine Angew. Math.} {\bf 690} (2014), 51--63. 

\bibitem{KnuthV4A}
D. E. Knuth, \newblock\newblock {\em The art of computer programming}, Vol. 4A, Addison-Wesley, 2011.

\bibitem{Krattenthaler}
C. Krattenthaler, \newblock Lattice path enumeration, \newblock in M. Bona ed., {\em Handbook of enumerative combinatorics}, CRC Press, 2015, 589--678.

\bibitem{Lindstrom1970}
B. Lindstr\"om, \newblock Conjecture on a theorem similar to Sperner's, \newblock in R. Guy ed., {\em Combinatorial structures and their applications}, Gordon and Breach, New York, 1970, 241.

\bibitem{MacDonald}
I. G. Macdonald, \newblock\newblock {\em Symmetric functions and Hall polynomials}, 2nd edition, Oxford University Press, 1995.

\bibitem{MSTT2018}
K. Manes, I. Tasoulas, A. Sapounakis, and P. Tsikouras, \newblock Counting pairs of noncrossing binary paths: A bijective approach, \newblock {\em Discrete Math.} {\bf 342}(2) (2019), 352--359.

\bibitem{MerrisRoby2005}
R. Merris and T. Roby, \newblock The lattice of threshold graphs, {\em JIPAM. J. Inequal. Pure Appl. Math.}, \newblock {\bf 6}(1) (2005), Article 2.

\bibitem{Pechenik2014}
O. Pechenik, \newblock Cyclic sieving of increasing tableaux and small Schr\"oder paths, \newblock {\em J. Combin. Theory Ser. A} {\bf 125} (2014), 357--378.

\bibitem{PSV2018}
T. Pressey, A. Stokke, and T. Visentin, \newblock Increasing tableaux, Narayana numbers and an instance of the cyclic sieving phenomenon, \newblock {\em Ann. Comb.} {\bf 20} (3) (2016), 609--621.

\bibitem{STT2006}
A. Sapounakis, I. Tasoulas, and P. Tsikouras, 
\newblock On the dominance partial ordering of Dyck paths,
\newblock \textit{J. Integer Seq.} {\bf 9} (2006), \#06.2.5.

\bibitem{OEIS}
N. J. A. Sloane, \newblock\newblock The on-line encyclopedia of integer sequences, \url{https://oeis.org}

\bibitem{Stanley1991}
R. P. Stanley, \newblock Some applications of algebra to combinatorics, 
\newblock {\em Discrete Appl. Math.} {\bf 34} (1991), 241--277.

\bibitem{Stanley2011}
R. P. Stanley, \newblock\newblock {\em Enumerative combinatorics}, Vol. 1, 2nd edition, Cambridge University Press, 2011.

\bibitem{MSTT2019}
I. Tasoulas, K. Manes, A. Sapounakis, and P. Tsikouras, \newblock
Chains with small intervals in the lattice of binary paths, \newblock \url{https://arxiv.org/abs/1911.10883}.

\bibitem{ThomasYong2009}
H. Thomas, A. Yong, \newblock A jeu de taquin theory for increasing tableaux, with applications to K-theoretic Schubert calculus, \newblock {\em Algebra Number Theory} {\bf 3}(2) (2009), 121--148.

\end{thebibliography}
\end{document}